\documentclass[smallextended]{svjour3}       
\usepackage{a4wide,graphicx}
\usepackage{amsmath, amsfonts, amssymb, mathtools}

\mathtoolsset{showonlyrefs}

\DeclareMathAlphabet{\mathbbold}{U}{bbold}{m}{n}

\usepackage{color,amsthm}

\usepackage[usenames,dvipsnames,svgnames,table]{xcolor}

\definecolor{bananamania}{rgb}{0.98, 0.91, 0.71}
\definecolor{bananayellow}{rgb}{1.0, 0.88, 0.21}
\newtheorem{thm}{Theorem}
\newtheorem{lem}[thm]{Lemma}
\newtheorem{prop}[thm]{Proposition}

\newtheorem{cor}[thm]{Corollary}

\theoremstyle{definition}
\newtheorem{defn}[thm]{Definition}

\newtheorem{rem}[thm]{Remark}
\newtheorem{exam}[thm]{Example}

\DeclareMathOperator{\spn}{span}

\DeclareMathOperator{\dist}{dist}

\DeclareMathOperator{\rint}{ri}
\DeclareMathOperator{\conv}{conv}

\DeclareMathOperator{\graph}{graph}
\DeclareMathOperator{\dis}{dis}

\let\leq\leqslant
\let\geq\geqslant
\let\emptyset\varnothing

\newcommand{\calH}{\ensuremath{\mathcal{H}}}

\newcommand{\calM}{\ensuremath{\mathcal{M}}}
\newcommand{\calN}{\ensuremath{\mathcal{N}}}

\newcommand{\calP}{\ensuremath{\mathcal{P}}}
\newcommand{\calQ}{\ensuremath{\mathcal{Q}}}
\newcommand{\calR}{\ensuremath{\mathcal{R}}}

\newcommand{\calW}{\ensuremath{\mathcal{W}}}

\newcommand{\barx}{\ensuremath{\bar{x}}}
\newcommand{\bary}{\ensuremath{\bar{y}}}
\newcommand{\barz}{\ensuremath{\bar{z}}}

\newcommand{\bmat}{\begin{matrix}}
\newcommand{\emat}{\end{matrix}}
\newcommand{\bbm}{\begin{bmatrix}}
\newcommand{\ebm}{\end{bmatrix}}
\newcommand{\bbma}{\begin{bmatrix*}[r]}
\newcommand{\ebma}{\end{bmatrix*}}
\newcommand{\bpm}{\begin{pmatrix}}
\newcommand{\epm}{\end{pmatrix}}
\newcommand{\bpma}{\begin{pmatrix*}[r]}
\newcommand{\epma}{\end{pmatrix*}}
\newcommand{\bvm}{\begin{vmatrix}}
\newcommand{\evm}{\end{vmatrix}}
\newcommand{\bse}{\begin{subequations}}
\newcommand{\ese}{\end{subequations}}
\newcommand{\beq}{\begin{equation}}
\newcommand{\eeq}{\end{equation}}
\newcommand{\ben}{\begin{enumerate}[1.]}

\newcommand{\een}{\end{enumerate}}

\newcommand{\beni}{\begin{enumerate}[(i)]}

\newcommand{\eeni}{\end{enumerate}}

\newcommand{\bena}{\begin{enumerate}[(a)]}

\newcommand{\eena}{\end{enumerate}}

\newcommand{\bit}{\begin{itemize}}
\newcommand{\eit}{\end{itemize}}
\newcommand{\bthe}{\begin{theorem}}
\newcommand{\ethe}{\end{theorem}}
\newcommand{\blem}{\begin{lemma}}
\newcommand{\elem}{\end{lemma}}
\newcommand{\bprop}{\begin{proposition}}
\newcommand{\eprop}{\end{proposition}}
\newcommand{\bex}{\begin{example}}
\newcommand{\eex}{\end{example}}
\newcommand{\bas}{\begin{assumption}}
\newcommand{\eas}{\end{assumption}}
\newcommand{\bre}{\begin{remark}}
\newcommand{\ere}{\end{remark}}
\newcommand{\bcor}{\begin{corollary}}
\newcommand{\ecor}{\end{corollary}}
\newcommand{\bdfn}{\begin{definition}}
\newcommand{\edfn}{\end{definition}}
\newcommand{\bcon}{\begin{conjecture}}
\newcommand{\econ}{\end{conjecture}}

\newcommand{\half}{\ensuremath{\frac{1}{2}}}
\newcommand{\inv}{\ensuremath{^{-1}}}
\newcommand{\nonempty}{\ensuremath{\neq\emptyset}}
\newcommand{\pset}[1]{\ensuremath{\{#1\}}}
\newcommand{\nset}[1]{\ensuremath{\{1,2,\ldots,#1\}}}

\newcommand{\set}[2]{\ensuremath{\{#1: #2\}}}

\newcommand{\abs}[1]{\ensuremath{| #1 |}}

\newcommand{\inn}[2]{\ensuremath{\left\langle #1 , #2 \right\rangle}}

\newcommand{\Lone}{\ensuremath{L_1}}

\newcommand{\Ltwo}{\ensuremath{L_2}}

\newcommand{\Linf}{\ensuremath{L_\infty}}

\newcommand{\abcd}{\ensuremath{(A,B,C,D)}}

\newcommand{\BP}{\noindent{\bf Proof. }}
\newcommand{\EP}{\hspace*{\fill} $\blacksquare$\bigskip\noindent}
\newcommand{\qand}{\quad\text{and}\quad}

\newcommand{\B}{\mathbb{B}}

\newcommand{\N}{\mathbb{N}}

\newcommand{\R}{\mathbb{R}}

\newcommand{\cF}{\mathcal{F}}

\newcommand{\cM}{\mathcal{M}}
\newcommand{\cN}{\mathcal{N}}

\DeclareMathOperator{\argmin}{argmin}

\DeclareMathOperator{\cl}{cl}

\DeclareMathOperator{\dom}{dom}
\DeclareMathOperator{\proj}{proj}

\DeclareMathOperator{\haus}{H}
\DeclareMathOperator{\im}{im}

\newcommand{\ul}{\underline}
\newcommand{\ol}{\overline}
\newcommand{\eps}{\varepsilon}
\newcommand{\setdef}[2]{\left\{\ #1\ \left|\ \vphantom{#1} #2\ \right.\right\}}

        {\hspace*{\fill}$\Box$\par}

\usepackage[prependcaption,colorinlistoftodos]{todonotes}

\usepackage{enumitem}

\begin{document}

%
%

\title{Convergence of proximal solutions for evolution inclusions with time-dependent maximal monotone operators
\thanks{The work of A. Tanwani is partially supported by ANR JCJC project ConVan with grant number ANR-17-CE40-0019-01.}
}
%
%
\author{Kanat Camlibel \and Luigi Iannelli \and  Aneel Tanwani}
\authorrunning{M.K.~Camlibel, L.~Iannelli and A. Tanwani} 
\institute{M.K.~Camlibel \at
              Bernoulli Institute of Mathematics, Computer Science, and Artificial Intelligence\\
              University of Groningen, The Netherlands\\
              \email{m.k.camlibel@rug.nl}%
           \and
           L.~Iannelli \at
           Department of Engineering\\
             University of Sannio in Benevento, Italy\\
             \email{luigi.iannelli@unisannio.it}
            \and
            A.~Tanwani \at
            Laboratory for Architecture and Analysis of Systems (LAAS), CNRS\\
            University of Toulouse, France\\
             \email{aneel.tanwani@laas.fr}%
}
\date{Received: date / Accepted: date}

\maketitle

\begin{abstract}
This article studies the solutions of time-dependent differential inclusions which is motivated by their utility in optimization algorithms and the modeling of physical systems. The differential inclusion is described by a time-dependent set-valued mapping having the property that, for a given time instant, the set-valued mapping describes a maximal monotone operator. By successive application of a proximal operator, we construct a sequence of functions parameterized by the sampling time that corresponds to the discretization of the continuous-time system. Under certain mild assumptions on the regularity with respect to the time argument, and using appropriate tools from functional and variational analysis, this sequence is then shown to converge to the unique solution of the original differential inclusion. The result is applied to develop conditions for well-posedness of differential equations interconnected with nonsmooth time-dependent complementarity relations, using passivity of underlying dynamics (equivalently expressed in terms of linear matrix inequalities).
\end{abstract}

\section{Introduction}

The theory of monotone operators emerged as an important area of research within the field of nonlinear analysis in early 1960's \cite{Kacu60,Mint62,Zara60}. Since then, we have seen applications of such operators in various disciplines, which include, but are not limited to, optimization algorithms, dynamical systems, and partial differential equations.
Recent articles \cite{attouch2019,Comb18,RyuBoyd16} provide an overview of monotone operators appearing in optimization algorithms. The relevance of set-valued mappings in dynamical systems observed in \cite{Brez73,More77}, where the differential inclusions with maximal monotone operators are analyzed. Even in the systems of partial differential equations, the appearance of these operators brings tractability to proving existence of solutions \cite{Brow63,LeraLion65,Show93,Zeid90}. Applications of dynamical systems with maximal monotone operators range from modeling traffic equilibrium \cite{pang2012} to electronics \cite{addi11}. Relatively modern texts on analysis of monotone operators are \cite{BausComb17,Phel93,Simo08}.

This article is focused on studying maximal monotone operators in the context of mathematical models for dynamical systems, and the central object of our study is to investigate conditions for existence of solutions to the differential inclusion
\begin{equation}\label{eq:sysIntro}
\dot x \in - F(t,x), \quad x(0) \in \dom F(0,\cdot),
\end{equation}
where $F : [0,\infty) \times \R^n \rightrightarrows \R^n$ has the property that, for each $t \geq 0$, $F(t, \cdot)$ is a maximal monotone operator.
In studying this generic class of systems, we will refer to other types of nonsmooth dynamical systems which can be recast in the form \eqref{eq:sysIntro}. From a theoretical point of view, most of the earlier work had focused on differential inclusions with {\em static} maximal monotone operators, which is very elegantly collected in \cite{Brez73}, or see \cite{PeypSori10} for a recent overview on this subject. Common techniques used in analyzing such systems are either based on regularization, or discretization. For the former one, the so-called Yosida-Moreau approximations provide a single-valued Lipschitz function with a regularization parameter, and as this parameter converges to zero, it is shown that the corresponding solutions converge to the solution of the original differential inclusion. The discretization techniques rely on constructing piecewise constant interpolations of the sequence of points obtained from some discrete system with a sampling parameter. As the sampling parameter converges to zero, the corresponding sequence of solutions is shown to converge to the actual solution.

To the best of our knowledge, the first attempts for studying inclusion \eqref{eq:sysIntro} with time-dependent operators $F(t,x)$, but with the domain of $F(t,\cdot)$ stationary for each $t \geq 0$, were carried out in \cite{Kato67}. Since then, several works have appeared which tackle dynamical systems with time-dependent multi-valued monotone operators. When $F(t,x)$ is the subdifferential of a time-dependent, proper, lower semicontinuous, and convex function $\varphi_t(\cdot)$, that is, $F(t,x) = \partial \varphi_t(x)$, then $F(t,\cdot)$ is a maximal monotone operator. Such systems, involving time-dependent subdifferentials, have been particularly studied in \cite{Arseni-Benou99,kandilakis96,kartsatos84,Kuttler00,otani94,yamazaki05} under varying degrees of regularity on the system data. Imposing further structure on the operator $F(t,\cdot)$, if we take $F(t,x) = \partial \psi_{S(t)}(x)$, where $S:[0,\infty) \rightrightarrows \R^n$ is closed and convex-valued mapping and $\psi$ is the indicator function associated with $S(t)$, then the resulting dynamics have been more commonly studied under the topic of {\em sweeping processes}. Starting from the seminal work of \cite{More77}, the research in this area has grown to study several generalizations of the fundamental model, see for example, the monographs \cite{Adly18,KunzMont00,Mont93,SiddManc02} for an overview, and the articles \cite{AdlyHadd18,AdlyHadd14,EdmoThib05,EdmoThib06,JourVilc18,KrejRoch11,Recu15} for more recent and focused expositions. Besides the cases where $F$ is expressed as a subdifferential of a convex function, certain classes of evolution variational inequalities \cite{BrogGoel11,PangStew08,TanwBrog18} can also be embedded in the framework of \eqref{eq:sysIntro}. 

While all these aforementioned works can be represented by \eqref{eq:sysIntro}, they also rely on the particular structure of the set-valued mapping in their problem description for analysis of existence of solutions. Notable exceptions in the literature, which address directly the system class \eqref{eq:sysIntro} are \cite{KunzMont97,Vlad91}. However, the regularity assumptions imposed in these works restrict the applicability of their results. Consequently, when applications of these dynamical models are studied, for example in control \cite{Briceno-Arias15,TanwBrog18}, the results that build on the works of \cite{KunzMont97,Vlad91} suffer similar limitations. Based on these observations, the motivation to study new set of conditions for existence of solutions to systems class~\eqref{eq:sysIntro}  arises and our aim in this paper is to provide mild (read as mildest possible) conditions on regularity of system data, which allow us to cover a possibly larger class of systems. Moreover, we can recover most (if not all) of the results on time-dependent and static case with our approach.




Our approach builds on using the time-stepping algorithm pioneered in \cite{More77}, which was also used for studying existence of solutions for system~\eqref{eq:sysIntro} in \cite{KunzMont97}. This algorithm constructs a sequence of solutions, where each element of the sequence is an interpolation of points obtained by applying the proximal operator associated with domain of the multivalued function appearing in \eqref{eq:sysIntro}. With the help of an academic counterexample, we show how the assumptions imposed in \cite{KunzMont97} fail to hold for a dynamical system described by time-dependent complementarity relations. We study existence of solutions under conditions which overcome such restrictions. The basic idea is to construct a sequence of solutions. To construct an element of this approximate  solution with a fixed sampling time, we first compute a set of points at sampled time instants by projecting the value of a certain function on the domain of the set-valued mapping. Using a novel interpolation technique among these discrete points, we obtain a sequence of absolutely continuous functions. Using the arguments based on Ascol\`a-Arzeli theorem, this sequence is shown to converge to an absolutely continuous function, which is then shown to be the unique solution of the original system. We generalize our result to the case where the right-hand side of \eqref{eq:sysIntro} has a single-valued Lipschitz vector field in addition to the set-valued maximal monotone operator.

Moreover, because of the relaxed nature of assumptions, our results provide a constructive framework for studying differential equations with complementarity relations. Such nonsmooth relations form a particular subclass of maximal monotone operators, and have been useful in modeling systems with piecewise affine characteristics \cite{j06,Brog03,BrogThib10}. Earlier work on complementarity systems has focused on linear dynamics coupled with static complementarity relations \cite{CamlHeem02,Heem00,schumacher2004,j16}. Lately, it was shown in \cite{CamlSchu16} that an interconnection of static complementarity relation with ordinary differential equations yields a differential inclusion with static maximal monotone operator. However, time-dependence in complementarity relations has not been easy to treat with existing frameworks. Inspired by the result in \cite{CamlSchu16}, we provide conditions under which it is possible to recast the interconnection of an ordinary differential equation with time-dependent complementarity relation in the form of a differential inclusion with time-dependent maximal monotone operator, for which the existence of solutions is being studied in this article.

The remainder of the article is organized as follows: In Section~\ref{sec:prelim}, we provide appropriate background material from set-valued and functional analysis. In Section~\ref{sec:maxMonTV}, a motivating example is provided to show how the current literature on differential inclusions with maximal monotone operators is inadequate for certain system classes. The main assumptions and the result is given in Section~\ref{sec:mainThm}, followed by a detailed proof in Section~\ref{sec:proof}. Section~\ref{s:ext} deals with extensions of the main existence/uniqueness result towards non-autonomous case as well as Lipschitzian perturbations. The results are then studied in the context of linear ordinary differential equations coupled with time-dependent maximal monotone relations in Section~\ref{sec:compSys}. Finally, the paper closes with some concluding remarks in Section~\ref{s:conc}.

\section{Preliminaries}\label{sec:prelim}
In this section, we introduce notational conventions that will be in force throughout the paper as well as auxiliary results that will be employed later.

\subsection{Vectors and matrices}
We denote the set of real numbers by $\R$, nonnegative real numbers by $\R_+$, $n$-vectors of real numbers by $\R^n$, and $n\times m$ real-valued matrices by $\R^{n\times m}$. 

To denote the scalar product of two vectors $x$, $y\in\R^n$, we use the notation $\inn{x}{y}:=x^\top  y$ where $x^\top $ denotes the transpose of $x$. The Euclidean norm of a vector $x$ is denoted by $\abs{x}:=\inn{x}{x}^\half$. For a subspace of $\calW$ of $\R^n$, $\calW^\perp$ denotes the orthogonal subspace, that is $\set{y\in\R^n}{\inn{x}{y}=0\text{ for all }x\in \calW}$.

We say that a (not necessarily symmetric) matrix $M\in\R^{n\times n}$ is {\em positive semi-definite\/} if $x^\top  Mx\geq 0$ for all $x\in\R^n$. We sometimes write $M\geq 0$ meaning that $M$ is positive semi-definite. Also, we say that $M$ is positive definite if $M> 0$ for all nonzero $x\in\R^n$.

\subsection{Convex sets and related notions}
The \emph{distance} of a point $x$ to a set $S$ is defined by $\dist(x,S)=\inf\set{\abs{x-y}}{y\in S}$. If the set $S$ is closed and convex then for each $x\in\R^n$ there exists a unique point $y\in S$ such that $\abs{x-y}=\dist(x,S)$. Such a point is called the \emph{projection} of $x$ onto the set $S$ and will be denoted by $\proj(x,S)$.

The \emph{Hausdorff distance} between two nonempty subsets of $\R^m$, say $S_1$ and $S_2$, is defined by:
\begin{align*}
d_{\haus}(S_1,S_2) := \max{\big\{\sup_{z_1 \in S_1} \dist(z_1,S_2), \sup_{z_2 \in S_2} \dist(z_2,S_1) \big\}}.
 \end{align*}
Since $\dist(x,S)=\dist\big(x,\cl(S)\big)$ for any point $x$ and nonempty set $S$, the Hausdorff distance is invariant under closure, that is
$$
d_{\haus}(S_1,S_2)=d_{\haus}\big(\cl(S_1),\cl(S_2)\big).
$$

In addition, if $y=\proj(x,\cl(S_2))$ for some point $x\in\cl(S_1)$, then we have
\beq\label{e:hausdorff conseq}
\abs{x-y}\leq \sup_{z \in S_1} \dist(z,S_2) \le d_{\haus}(S_1,S_2).
\eeq

\subsection{Set-valued mappings}

Let $F:\R^m\rightrightarrows\R^n$ be a set-valued mapping, that is $F(x)\subseteq\R^n$ for each $x\in\R^m$. We define its domain, image, and graph, respectively, as follows:
\begin{gather*}
\dom F=\set{x}{F(x)\neq\emptyset}\\
\im F=\set{y}{\text{there exists }x\text{ such that }y\in F(x)}\\
\graph F =\set{(x,y)}{y\in F(x)}.
\end{gather*}
The inverse mapping $F\inv:\R^n\rightrightarrows\R^m$ is defined by $F\inv(y)=\set{x}{y\in F(x)}$.

In what follows we introduce a certain notion of continuity for set-valued mappings of a real variable. For a more detailed/general treatment we refer to \cite[Chp. 4 and 5]{Rockafellar98}.

Let $\calN_\infty^\#$ denote the set of all subsequences of $\N$. For a sequence of sets $(S_\ell)_{\ell\in \N}$ in $\R^q$, the {\em outer limit\/} is defined as the set
$$
\limsup_{\ell\rightarrow\infty} S_\ell:=\setdef{\xi}{\exists\,N\in\calN_\infty^\#\text{ and }\xi_\ell\in S_\ell\,\,\forall\,\ell\in N, \text{ s.t. } \xi_\ell\overset{N}{\rightarrow}\xi}.
$$
For a given set-valued mapping $G:[0,T]\rightrightarrows\R^q$ for some $T > 0$, we define
$$
\limsup_{t\searrow t^*}G(t):=\bigcup_{t_\ell\searrow t^*}\limsup_{\ell\rightarrow\infty} G(t_\ell).
$$
It is known from \cite[p. 152]{Rockafellar98} that
$$
\limsup_{t\searrow t^*}G(t) =\setdef{y\in\R^q}{\begin{aligned} & \exists\,(t_\ell,y_\ell)_{\ell\in\N}\subset[0,T]\times\R^q \text{ satisfying } y_\ell\in G(t_\ell), \\ & {t_{\ell} \ge t}, \text{ and }\lim_{\ell\rightarrow\infty}(t_\ell,y_\ell)=(t^*,y)\end{aligned}}.
$$
{We say that $G$ is outer semicontinuous from right at $t^*\in[0,T]$ if
$$
\limsup_{t\searrow t^*}G(t)\subseteq G(t^*).
$$
In case $G$ is outer semicontinuous from right at every $t^*\in[0,T]$, we say that $G$ is outer semicontinuous from right on $[0,T]$.
}

\subsection{Maximal monotone operators}
Throughout the paper, we are interested in maximal monotone set-valued mappings. A set valued-mapping $F:\R^n\rightrightarrows\R^n$ is said to be {\em monotone\/} if
\beq
\inn{x_1-x_2}{y_1-y_2}\geq 0
\eeq
for all $(x_i,y_i)\in\graph(F)$. It is said to be {\em maximal monotone\/} if no enlargement of its graph is possible in $\R^n\times\R^n$ without destroying monotonicity. We refer to \cite{Brez73} and \cite{Rockafellar98} for detailed treatment of maximal monotone mappings.

If $F$ is maximal monotone, then it is closed and convex-valued, that is, $F(x)$ is a closed convex set for all $x\in\dom(F)$. This enables us to define
the \emph{minimal section} of a maximal monotone mapping $F$ by
$$
F^0(x):=\proj(0,F(x))
$$
for $x\in\dom(F)$. Clearly, $F^0(x)$ is the least-norm element of the closed convex set $F(x)$, that is $\abs{F^0(x)}\leq \abs{y}$ for all $y\in F(x)$.

The \emph{resolvent } $J_\lambda$ and \emph{Yosida approximation} $Y_\lambda$ of $F$ are defined by
$$
J_\lambda:=(I +\lambda F)^{-1}
\qand
Y_\lambda:=\frac{1}{\lambda}(I -J_\lambda)
$$
for $\lambda >0$ where $I$ denotes the identity operator.

The following proposition collects some well-known facts (see e.g. \cite{Brez73}) that will be employed in the sequel.

\begin{prop}\label{prop:maxmon-resolvent-Yosida}
Suppose that $F:\R^n\rightrightarrows\R^n$ is a maximal monotone set-valued mapping. Then, the following statements hold for all $\lambda>0$:
\begin{enumerate}[label=\roman*.]
\item $\dom J_\lambda =\R^n$.
\item $J_\lambda$ is single-valued and non-expansive, that is $\abs{J_\lambda (x_1)- J_\lambda (x_2)} \leq \abs{x_1-x_2}$ for all $x_1$, $x_2\in\R^n$.
\item $\lim_{\lambda\rightarrow 0}J_\lambda(x)=x$ for all {$x\in\dom(F)$}.
\item $Y_\lambda$ is maximal monotone and $\lambda\inv-$Lipschitzian.
\item $Y_\lambda (x)\in F\big(J_\lambda(x)\big)$ for all $x\in\R^n$.
\item For all $x\in \dom F $, $\abs{Y_\lambda (x)}$ is nonincreasing in $\lambda$, $\lim_{\lambda\rightarrow 0}\abs{Y_\lambda (x)}=\abs{F^0 (x)}$, and $\abs{Y_\lambda (x)}\leq \abs{F^0 (x)}$.
\end{enumerate}
\end{prop}

Given two maximal monotone mappings, the {\em pseudo-distance\/} between them, introduced in \cite{Vlad91}, is defined as follows:

\begin{defn}\label{def:dis}
	The {\em pseudo-distance\/} between two maximal monotone mappings $F_1$ and $F_2$ is defined by
$$
	\dis(F_1,F_2):=\sup_{
		\small \begin{array}{l}
			x_1 \in \dom(F_1), y_1 \in F_1(x_1)\\ x_2 \in \dom(F_2), y_2 \in F_2(x_2)
		\end{array}}{\dfrac{\langle y_1-y_2, x_2-x_1\rangle}{1+\abs{y_1}+\abs{ y_2}}}.
$$
\end{defn}

The following lemma relates the Hausdorff distance between the domains of two maximal monotone operators with their pseudo-distance.

\begin{lem}[\cite{Vlad91}]
For any pair of maximal monotone mappings $F_1$ and $F_2$, it holds that
$$d_{\haus}(\dom(F_1),\dom(F_2))\leq \dis(F_1,F_2).$$
\end{lem}

Based on the pseduo-distance defined in Definition~\ref{def:dis}, one can introduce a notion of continuity for time-dependent maximal monotone operators as follows.

\begin{defn}[Absolute continuity,~\cite{Vlad91}]\label{d:vlad}
Let $F:[0,T]\times\R^n\rightrightarrows\R^n$ be a time-dependent set-valued mapping such that $F(t,\cdot)$ is maximal monotone for each $t \in [0,\,T]$. We say that $t\mapsto F(t,\cdot)$ is {\em absolutely continuous\/} on $[0,T]$ if there exists a nondecreasing absolutely continuous function $\varphi\; : \; [0,\,T] \rightarrow \R$ such that
	\begin{align*}
		\dis\big(F(t,\cdot),F(s,\cdot)\big)\leq \varphi(t)-\varphi(s) \quad \forall\,s,t\text{ with } 0\leq s \leq t\leq T.
	\end{align*}
\end{defn}

\subsection{Function spaces}
The set of absolutely continuous, integrable, and square integrable functions defined from the interval $[t_1,t_2]$ with $t_1<t_2$ to $\R^n$ are denoted, respectively, by $AC([t_1,t_2],\R^n)$, $\Lone([t_1,t_2],\R^n)$, and $\Ltwo([t_1,t_2],\R^n)$. Unless specified otherwise, we use the term {\em almost everywhere} with respect to Lebesgue measure, that is, a property holds almost everywhere on a set $X \subset \R^n$, if it holds on every subset of $X$ with nonzero Lebesgue measure.

Convergence of family of functions will play an important role in the sequel. For the sake of completeness, we state the well-known (see e.g. \cite{rudin}) Arzel\'{a}-Ascoli theorem for which we need some nomenclature.

Consider a collection $\cF$ of functions $f:[0,T] \to \R^n$. We say that $\cF$ is equicontinuous if for every $\eps > 0$, there exists a $\delta > 0$ such that $\vert f(t) - f(s) \vert < \eps$ for every $f \in \cF$ and each $s,t$ satisfying $\vert t - s \vert < \delta$. We say that $\cF$ is pointwise bounded if for every $t \in [0,T]$, there exists an $M_t < \infty$ such that $\vert f(t) \vert \le M_t$ for every $f \in \cF$.

\begin{thm}[Arzel\'{a}-Ascoli]\label{thm:ascArz}
Suppose that $\cF$ is pointwise bounded equicontinuous collection of functions $f:[0,T] \to \R^n$. Every sequence $\{f_n\}$ in $\cF$ has a subsequence that converges uniformly on every compact subset of $[0,T]$.
\end{thm}

The following elementary results will be used later.

\begin{lem}\label{l:der lim}
Let $x:[0,T]\rightarrow\R^n$ be a function and $t^*\in(0,T]$ be such that $\dot{x}(t^*)$ exists. Suppose that $\pset{t_k}$ and $\pset{\tau_k}$ are two sequences such that $0\leq t_k\leq t^*\leq \tau_k\leq T$ and $t_k<\tau_k$ for all $k$ and $\lim_{k\uparrow\infty} t_k=\lim_{k\uparrow\infty}\tau_k=t^*$. Then, the sequence $\frac{x(\tau_k)-x(t_k)}{\tau_k-t_k}$ convergences to $\dot{x}(t^*)$ on a subsequence.
\end{lem}

\BP
Observe that
$$
\frac{x(\tau_k)-x(t_k)}{\tau_k-t_k}=\frac{x(\tau_k)-x(t^*)}{\tau_k-t^*}\frac{\tau_k-t^*}{\tau_k-t_k}+\frac{x(t^*)-x(t_k)}{t^*-t_k}\frac{t^*-t_k}{\tau_k-t_k}.
$$
Since $0\leq \frac{\tau_k-t^*}{\tau_k-t_k}\leq 1$ and $0\leq \frac{t^*-t_k}{\tau_k-t_k}\leq 1$, both must converge on a (common) subsequence. The rest follows from the hypothesis that $\dot x(t^*)$ exists.
\EP

\begin{lem}\label{lem:xderpsiy}
Suppose that a sequence of functions $(y_\ell)_{\ell\in\N}$ weakly converges to $y$ in $\Ltwo(d\psi,[0,T],\R)$, for some $\psi \in AC([0,T], \R)$. Let $(x_\ell)_{\ell\in\N}$ be a sequence of absolutely continuous functions such that it converges uniformly to $x \in AC([0,T],\R^n)$, and $\dot x_\ell(t) =\dot \psi(t) y_\ell(t)$, for each $t \in \Gamma:= \bigl\{t\in [0,T] \, \vert \, x_\ell, x$ and $\psi$ are differentiable at  $t\bigr\}$. Then, it holds that $\dot x (t) = \dot \psi (t) y(t)$ for almost every $t\in \Gamma$.
\end{lem}

\BP
Define the function $\xi:[0,T]\rightarrow\R^n$ by 
\beq\label{eq:xiLim}
\xi(t)=x(0)+\int_0^t y(s)\dot{\psi}(s)\,ds
\eeq
for $t\in[0,T]$. For every $\eta\in\R^n$, we have
$$
\inn{\eta}{x_{\ell}(t)}=\inn{\eta}{x_0}+\int_0^t\inn{\eta}{y_{\ell}(s)}\,\dot \psi (s)\,ds
$$
for all $\ell\in \N$ and
$$
\inn{\eta}{\xi(t)}=\inn{\eta}{x_0}+\int_0^t\inn{\eta}{y(s)}\,\dot{\psi}(s)\,ds.
$$
Since $(y_{\ell})_{\ell\in \N}$ weakly converges to $y$, we have that $\big(\inn{\eta}{x_{\ell}(t)}\big)_{\ell\in \N}$ converges to $\inn{\eta}{\xi(t)}$ for every $t\in[0,T]$ and every $\eta\in\R^n$. This means that $\big(x_{\ell}(t)\big)_{\ell\in \N}$ converges to $\xi(t)$ for every $t\in[0,T]$. Hence, we see that $\xi(t)=x(t)$ for all $t\in[0,T]$ since $\big(x_{\ell}\big)_{\ell\in \N}$ uniformly converges to $x$. Therefore, \eqref{eq:xiLim} yields
\beq\label{e:from xi to x}
x(t)=x_0+\int_0^t y(s)\dot{\psi}(s)\,ds.
\eeq
In other words, 
\beq
\dot{x}(t)=\dot{\psi}(t)y(t)
\eeq
for almost every $t\in\Gamma$.
\EP

For the next two statements, we recall that two measures are absolutely continuously equivalent if each one is absolutely continuous with respect to the other one.

\begin{lem}\label{l:values of weak limit}
Let $f_\ell:[0,T]\rightarrow\R$ be a sequence of functions with $\ell\in\N$ such that $\abs{f_\ell(t)}\leq 1$ for all $\ell\in\N$ and $t\in[0,T]$. Suppose that the sequence $(f_\ell)_{\ell\in\N}$ weakly converges to $f$ in $\Ltwo(d\mu,[0,T],\R)$ where $d\mu$ is absolutely continuously equivalent to Lebesgue measure. Then, 
$$
f(t)\in[\smashoperator{\liminf_{\ell{\rightarrow}\infty}}f_\ell(t),\smashoperator{\limsup_{\ell{\rightarrow}\infty}}f_\ell(t)]
$$
for almost all $t\in[0,T]$.
\end{lem}
\BP
Let $k\geq 1$ and define $g^k_\ell(t):=\sup_{q\geq k}f_q(t)-f_{\ell+k}(t)$. Note that $(g^k_\ell)_{{\ell}\in\N}$ weakly converges in $\Ltwo(d\mu,[0,T],\R)$ to $g^k$ given by $g^k(t):=\sup_{q\geq k}f_q(t)-f(t)$. Since $g^k_\ell$ is nonnegative for all $\ell\in\N$ and $t\in[0,T]$, $g^k$ must be nonnegative for almost all $t\in[0,T]$. This means that $f(t)\leq \sup_{q\geq k}f_q(t)$ for almost all $t\in[0,T]$. Hence, $f(t)\leq \limsup_{\ell{\rightarrow}\infty}f_\ell(t)$ for almost all $t\in[0,T]$. Applying the same arguments to the sequence $(-f_\ell)_{\ell\in\N}$, we can obtain $f(t)\geq \liminf_{\ell{\rightarrow}\infty}f_\ell(t)$ for almost all $t\in[0,T]$. 
\EP

\begin{lem}\label{l:values of weak limit-general}
Let $y_\ell:[0,T]\rightarrow\R^q$ be a sequence of functions with $\ell\in\N$ such that $\abs{y_\ell(t)}\leq 1$ for all $\ell\in\N$ and $t\in[0,T]$. Also let $\big(S_\ell(t)\big)_{\ell\in\N}$ be a sequence of sets in $\R^q$ with $\ell\in\N$ and $t\in[0,T]$ such that $y_\ell(t)\in S_\ell(t)$ for all $\ell\in\N$ and $t\in[0,T]$. Suppose that $(y_\ell)_{\ell\in\N}$ weakly converges to $y$ in $\Ltwo(d\mu,[0,T],\R^q)$ where $d\mu$ is absolutely continuously equivalent to Lebesgue measure. Then, 
$$y(t)\in\cl\Big(\conv\big(\limsup_{\ell\rightarrow\infty}S_\ell(t)\big)\Big)$$
for almost all $t\in[0,T]$.
\end{lem}

\BP
Let $S(t)=\cl\Big(\conv\big(\limsup_{\ell\rightarrow\infty}S_\ell(t)\big)\Big)$ for $t\in[0,T]$. It follows from \cite[Cor. 4.11]{Rockafellar98} that $S(t)\nonempty$ for each $t\in[0,T]$. Let $\Gamma=\set{t\in[0,T]}{y(t)\notin S(t)}$. Define the function $z:[0,T]\rightarrow \R$ by
$$
z(t)=\proj\big(y(t),S(t)\big)
$$
for all $t\in \Gamma$ and $z(t)=0$ for all $t\in[0,T]\setminus \Gamma$. Note that $\abs{y(t)-z(t)}>0$ for all $t\in \Gamma$. Also, we have $z\in\Linf([0,T],\R^n)$ since $S(t)$ contains an element in the unit ball of $\R^n$ for all $t\in[0,T]$. Now, define functions $a:[0,T]\rightarrow \R^n$ and $b:[0,T]\rightarrow \R$ by
$$
a(t)=\frac{y(t)-z(t)}{\abs{y(t)-z(t)}}\qand b(t)=\inn{\frac{y(t)-z(t)}{\abs{y(t)-z(t)}}}{\frac{y(t)+z(t)}{2}}
$$
for all $t\in \Gamma$ and $a(t)=0$, $b(t)=0$ for all $t\in[0,T]\setminus \Gamma$. For all $t\in \Gamma$, the hyperplane $\calH_t=\set{\eta}{\inn{a(t)}{\eta}=b(t)}$
strictly separates the set $S(t)$ and the point $y(t)$, that is 
\beq\label{e:strict sep}
\inn{a(t)}{y(t)}<b(t)<\inn{a(t)}{z}
\eeq 
for all $z\in S(t)$ (see e.g. \cite[Prop. 1.5.3]{bertsekas}).
Note that $a\in\Linf([0,T],\R^n)$, and since $y \in \Ltwo(d\mu,[0,T],\R^n)$ and $z \in \Linf([0,T],\R^n)$, it follows that $b\in\Ltwo(d\mu,[0,T],\R)$.
 Therefore, the function $t\mapsto \inn{a(t)}{w(t)}$ belongs to $\Ltwo(d\mu,[0,T],\R)$ for every $w\in\Ltwo(d\mu,[0,T],\R^n)$. For each $\ell\in \N$, define $\zeta_\ell:[0,T]\rightarrow\R$ with $\zeta_\ell(t)=\inn{a(t)}{y_{\ell}(t)}$ for all $t\in \Gamma$ and $\zeta_\ell(t)=0$ for all $t\in[0,T]\setminus \Gamma$. Then, we see that $(\zeta_\ell)_{\ell\in \N}$ weakly converges to $\zeta$ given by $\zeta(t)=\inn{a(t)}{y(t)}$ for all $t\in \Gamma$ and $\zeta(t)=0$ for all $t\in[0,T]\setminus \Gamma$. From Lemma~\ref{l:values of weak limit}, we see that
\beq\label{e:zetas}
\zeta(t)\in[\smashoperator{\liminf_{\ell{\rightarrow}\infty}}\zeta_\ell(t),\smashoperator{\limsup_{\ell{\rightarrow}\infty}}\zeta_\ell(t)]
\eeq
for almost all $t\in[0,T]$. Since limit inferior (superior) can be obtained as the limit of a subsequence, \eqref{e:zetas} implies that for almost all $t\in \Gamma$
$$
\inn{a(t)}{y(t)}\in[\inn{a(t)}{\ul y(t)},\inn{a(t)}{\ol y(t)}]
$$
where $\ul y(t)$ and $\ol y(t)$ belong to $S(t)$. Together with the second inequality in \eqref{e:strict sep}, this yields
\beq\label{e:towards the end}
\inn{a(t)}{y(t)}>b(t)
\eeq
for almost all $t\in \Gamma$. In view of the first inequality in \eqref{e:strict sep}, this means that $\Gamma$ is a zero measure set. As such, we can conclude that 
$$y(t)\in S(t)$$
for almost all $t\in[0,T]$.
\EP

\section{Differential inclusions with maximal monotone mappings}\label{sec:maxMonTV}
Our goal is to study the existence of solutions to the differential inclusion
\begin{equation}\label{eq:incGen}
\dot{x}(t)\in -F\big(t,x(t)\big),\quad x(0)=x_0, \quad {t \in [0,T],} 
\end{equation}
where $F(t,\cdot):\R^n \rightrightarrows \R^n$ is maximal monotone for all {$t\in [0,T]$}. We say that $x\in AC([0,T],\R^n)$ is a solution of \eqref{eq:incGen} if $x(t) \in \dom F(t,\cdot)$ and $x$ satisfies \eqref{eq:incGen} for almost all $t\in [0,T]$.

%
%
%

\subsection{Related frameworks and their limitations}\label{sec:connections}
Historically, the evolution inclusions given in \eqref{eq:incGen} have been a subject of research in mathematical community in different eras. However, the solutions to such equations have been proposed under rather strict conditions. Here, we provide a brief of list of the main results that exist concerning the existence and uniqueness of solutions for such systems. {We also refer the reader to \cite[Chapter II]{HuPapa00} for an overview of the special case when $F(t,\cdot)$ is the subdifferential of a convex function.}\\

\noindent{\em Single-valued operators with fixed domain \cite{Kato67}:} The earliest results on solutions of dynamical systems \eqref{eq:incGen} with time-dependent maximal monotone relations were proposed in \cite{Kato67}. The author focused on the case where $F(t,\cdot):\R^n \to \R^n$ is single-valued and $F(\cdot,x)$ is Lipschitz continuous, uniformly in $x$. The major restriction imposed here is that
\begin{equation}\label{eq:constDom}
\dom F(t,\cdot) = \dom F(0,\cdot), \quad \forall \, t \geq 0.
\end{equation}
Under these conditions, there exists a Lipschitz continuous $x:\R_+ \to \R^n$ such that \eqref{eq:incGen} holds for Lebesgue almost every $t \geq0 $, and $x(t) \in \dom F(t,\cdot)$ for each $t \geq 0$.\\

\noindent{\em Static maximal monotone operators (\`a la Br\'ezis) with additive inputs \cite{Brez73}:} In the classical book \cite{Brez73} dealing with differential inclusions with maximal monotone operators, we can find results dealing with inclusions of the form
\[
\dot x(t) \in -F(t,x) = -A\big(x(t)\big) + \gamma \, x(t) - u(t)
\]
where $A$ is maximal monotone, $\gamma > 0$ is a scalar, and $u:[0,\infty) \to \R^n$ is absolutely continuous. For such systems $\dom F(t,\cdot)  = \dom  A$ for each $t \geq 0$, that is, the domain of the multivalued operator is independent of time.\\

\noindent{\em Dissipative operators \cite{Pave87}:} Building up on the work of Kato \cite{Kato67} and Br\'ezis \cite{Brez73}, we find results on evolution equations built on convergence of certain discrete approximations in \cite{Pave87}. When the results appearing in this line of work are applied to the maximal monotone case given in \eqref{eq:incGen}, it turns out that such results also require the strong assumption \eqref{eq:constDom}, where the domain of the operator does not change with time \cite[Chap.~1, Sec.~4]{Pave87}.\\

\noindent{\em Moreau's sweeping process \cite{More77}:} The first real contribution in the literature with time-dependent domains is observed in the seminal work of \cite{More77}. The systems studied here  within the umbrella of sweeping processes comprise differential inclusions with a special conic structure. We introduce a set-valued mapping $S:\R_+ \rightrightarrows \R^n$, and let $\cN_{S(t)}(x)$ denote the normal cone to the set $S(t)$ at a point $x \in S(t)$. The proposed system class is then described as:
\begin{equation}\label{eq:sweepFirst}
\dot x \in -F(t,x) :=-\cN_{S(t)} (x), \qquad x(0) \in S(0).
\end{equation}
Thus, for each $t$ and $x$, $F(t,x)$ is a closed convex cone described by the subdifferential of the indicator function of $S(t)$, and hence $F(t,\cdot)$ is maximal monotone. Here, we see that $\dom(F(t,\cdot)) = S(t)$ and since $S$ is time-dependent, the domain is allowed to vary with time. To describe the regularity imposed on $F(\cdot,x)$ with respect to time, we consider the Hausdorff distance, and in the simplest instances, it is assumed that, for every $t_1,t_2 \geq 0$
\[
d_{\haus} \big(\dom F(t_1,\cdot), \dom F(t_2,\cdot)\big) =
d_{\haus} \big(S(t_1), S(t_2)\big) \le L \vert t_1 - t_2 \vert,
\]
that is the Hausdorff distance between the domains of $F(t,\cdot)$ is bounded by a Lipschitz continuous function of time. Under these assumptions, there exists a unique solution to \eqref{eq:sweepFirst} which is Lipschitz continuous.
Different variants of this framework were then derived depending on how the Hausdorff distance varies with time, or whether we can relax the convexity assumption on $S(t)$ while preserving some nice properties of the subdifferential of the indicator function for that set.
In short, sweeping processes provide the first instance in the literature on inclusions with a particular kind of maximal monotone operators which depend on time, and the corresponding domain may vary.\\

\noindent{\em Maximal monotone operators with time-dependent domain \cite{Vlad91,KunzMont97}:} As a generalization of the sweeping process, Vladimirov \cite{Vlad91} studied evolution inclusions where time-dependent domains were considered, with the hypothesis that the set-valued mapping $F(t,\cdot)$ is just maximal monotone for each $t \geq 0$, without any further structural or geometrical assumption.
However, a very strong regularity assumption was imposed with respect to the pseudo-distance given in Definition~\ref{def:dis}. In particular, the mapping $F(t,\cdot)$ is required to be uniformly continuous, that is, there exists a sequence of piecewise constant operators $F_i: [0,T] \times \R^n \rightrightarrows \R^n$ such that for each $t \in [0,T]$
\[
\lim_{i \to \infty} \dis \big(F_i (t, \cdot),F(t,\cdot)\big) = 0.
\]
Kunze and Monteiro-Marques \cite{KunzMont97} then generalized this line of work to consider systems where the regularity with respect to time can be relaxed, so that the pseudo-distance between $F(t_1,\cdot)$ and $F(t_2,\cdot)$ is bounded by $\vert \mu(t_1) - \mu(t_2) \vert$ for some function of bounded variation $\mu: [0,T] \to \R $. Certain results developed in the context of sweeping processes are thus covered within this framework.
The work started by Vladimirov, and later generalized to some extent by Kunze and Monteiro-Marques, indeed is an attempt to deal with differential inclusions with most general time-dependent maximal monotone operators. However, they impose very strong assumptions in deriving their results which make their applicability somewhat restrictive. Indeed, as we show in the next section, strong continuity assumption is not necessary. 

\subsection{Motivation}\label{ssec:motiv}

A primary motivation for looking at inclusions of type~\eqref{eq:incGen} comes from differential equations where certain variables are related by a maximal monotone operator. In particular, consider systems described by
\begin{subequations}\label{eq:sysMotFirst}
\begin{align}
\dot x(t) &= A x(t) + B z(t) \\
w(t) & = Cx(t) + Dz(t) + v(t)\\
&w(t) \in \cM\big(-z(t)\big)
\end{align}
\end{subequations}
where $x\in\R^n$, $(z,w)\in\R^m\times\R^m$, $v\in\R^m$, the matrices $(A,B,C,D)$ have appropriate dimensions and $\cM:\R^{m} \rightrightarrows \R^m$ is a maximal monotone operator.

Systems of the form \eqref{eq:sysMotFirst} can be alternatively described by \eqref{eq:incGen} where
\beq\label{e:F t x}
F (t,x):=-Ax + B (\cM + D)^{-1}\big(Cx + v(t)\big).
\eeq
By invoking \cite[Theorem~2]{CamlSchu16}, one can show that $F (t,\cdot)$ is maximal monotone for each $t\in\R_+$ under certain assumptions. Regularity with respect to time is critical here. In the works of \cite{Vlad91,KunzMont97}, existence and uniqueness of solutions to \eqref{eq:incGen} is established under the assumption of absolute continuity in the sense of Definition~\ref{d:vlad}. However, the mapping $t\mapsto F(t,\cdot)$ defined by \eqref{e:F t x} does not, in general, enjoy absolute continuity with respect to pseudo-distance even if $v$ is absolutely continuous. This is seen in the following example.
\begin{exam}\label{ex} Consider a system of the form \eqref{eq:sysMotFirst} where $n=1$, $m=2$, 
$$
	A=0,\
	B=C^\top =\begin{bmatrix}
	0 &\,& 1
	\end{bmatrix}, \
	D=\bbma
	0 &\,& 1\\
	-1 &\,& 0
	\ebma,
$$
and $\cM:\R^2\rightrightarrows\R^2$ is the set-valued mapping given by $\cM(\zeta)=\set{\eta}{\eta\geq 0,\, \zeta \leq 0,\text{ and }\inn{\eta}{\zeta}=0}$. By invoking \cite[Theorem~2]{CamlSchu16}, it can be verified that the corresponding set-valued mapping $F(t,\cdot)$ as defined in \eqref{e:F t x} is maximal monotone for each $t$.  Let $v:[0,T] \to \R^2$ be an absolutely continuous function such that for some $t_1, t_2 \in [0,T]$, we have
\[
	v(t_1) =\begin{bmatrix} 0 \\ 0 \end{bmatrix}\qand
	v(t_2) =\begin{bmatrix} -1\\ 0 \end{bmatrix}.
\]
Let $F_i:= F (t_i,\cdot)$ with $i = 1,2$. It can be verified that
$$
0\in F_1(\rho+1)\qand 1\in F_2(0)
$$
for any $\rho\geq 0$. From Definition~\ref{def:dis}, we get
$$\label{eq:expVlad}
   	\dis(F_1,F_2) = \sup_{\footnotesize\begin{aligned}x \in \dom(F_1), y \in F_1(x), \\ \xi \in \dom(F_2), \zeta \in F_2(\xi)\end{aligned}}{\dfrac{\langle \zeta-y, x-\xi\rangle}{1+\abs{y}+\abs{\zeta}}}\geq \frac{\rho+1}{2}.
	$$
Since the righthand side is not bounded, we can conclude that set-valued mapping $F(t,\cdot)$ is not absolutely continuous in the sense of Definition~\ref{d:vlad}. However, existence and uniqueness of solutions for this example would follow from our main results. Indeed, this example satisfies the hypothesis of Theorem~\ref{e: linear case-comp}. 
\end{exam}

\section{Main results}\label{sec:mainThm}
The main goal of this paper is to investigate conditions (weaker than those of \cite{Vlad91,KunzMont97}) that guarantee existence of solutions to \eqref{eq:incGen}.  The uniqueness of solutions for a fixed initial condition follows easily from the maximal monotone property of the right-hand side of \eqref{eq:incGen}.

\subsection{Existence and uniqueness of solutions to \eqref{eq:incGen}}
To state the main result of our paper, we introduce the following assumptions, where $T>0$ and $\mathbb{B}^n(r):=\set{x\in\R^n}{\abs{x}\leq r}$.

\begin{enumerate}[label=(A\arabic*),leftmargin=*]
\item\label{ass:basicF} For each $t \in [0,T]$, the operator $ F (t,\cdot)$ is maximal monotone.
{\item\label{ass:new} The set-valued mapping $t\mapsto \graph F(t,\cdot)$ is outer semicontinuous from right on $[0,T]$.}
\item\label{ass:dist} There exists a nondecreasing function $\varphi\in AC([0,T],\R)$ such that 
$$
 \sup_{z \in \dom F(s,\cdot)}\dist \big(z , \dom F (t,\cdot)\big)  \leq \varphi(t)- \varphi(s), \quad \forall s,t \text{ with } 0\leq s\leq t \leq T.
$$
\end{enumerate}%
{
\begin{enumerate}[label=(LG),leftmargin=*]
\item\label{ass:least} There exists $\sigma \in\Lone([0,T],\R_+)$ such that 
$$
\abs{F ^0(t,x)}\leq \sigma (t) (1+\abs{x})
$$
for all $x \in  \dom{F(t,\cdot)}$ with $t\in[0,T]$. 
\end{enumerate}
\begin{enumerate}[label=(YB),leftmargin=*]
\item\label{ass:static} There exist a nondecreasing function $\theta:[0,T]\rightarrow\R$, and a scalar $\Lambda > 0$, such that for every $0 < \lambda < \Lambda $,
$$
\abs{Y_\lambda(s+\lambda,x)} \leq \theta(s+\lambda)-\theta(s)+\abs{F^0(s,x)}
$$
for all $s \text{ with } 0\leq s \leq T - \lambda$, and $x\in \dom F(s,\cdot)$ where $Y_{\lambda}(t,\cdot)$ denotes the Yosida approximation of $F(t,\cdot)$.
\end{enumerate}
The result on existence and uniqueness of solutions now follows.

\begin{thm}\label{thm:main}
Consider the system \eqref{eq:incGen} and assume that \ref{ass:basicF}--\ref{ass:dist} hold. If the linear growth assumption \ref{ass:least} holds, then there exists a unique solution of \eqref{eq:incGen} for every $x_0 \in \cl\big(\dom  F (0,\cdot)\big)$. In case the Yosida approximation bound assumption \ref{ass:static} holds, then there exists a unique solution of \eqref{eq:incGen} for every $x_0 \in \dom  F (0,\cdot)$.
\end{thm}

\begin{rem} The statement \ref{ass:static} in the description of Theorem~\ref{thm:main} can be checked without the explicit knowledge of the least norm element $\abs{F^0(s,x)}$ for each $s \ge 0$. In particular, consider the following statement:
\begin{enumerate}[label=(YB'),leftmargin=*]
\item\label{ass:static2} There exist a nondecreasing function $\theta:[0,T]\rightarrow\R$, and $\Lambda > 0$,  such that for every $0 < \lambda < \Lambda $,
$$
\abs{Y_\lambda(s+\lambda,x)} \leq \theta(s+\lambda)-\theta(s)+\abs{Y_\lambda(s,x)}
$$
for all $s \text{ with } 0\leq s \leq T - \lambda$, and $x\in \dom F(s,\cdot)$.
\end{enumerate}
Then, because of Proposition~\ref{prop:maxmon-resolvent-Yosida}, item (vi), \ref{ass:static2} implies \ref{ass:static}.
\end{rem}

\begin{rem}
If we consider the static maximal monotone operator $F:\R^n \rightrightarrows \R^n$, and the differential inclusion $\dot x \in -F(x)$, then the classical result on existence of solutions of such differential inclusions, such as \cite[Theorem~3.1]{Brez73}, follows directly from Theorem~\ref{thm:main}. In that case, conditions \ref{ass:new} and \ref{ass:dist} hold trivially, because the operator is time-independent and one can check that \ref{ass:static} holds with a constant function $\theta$, because of Proposition~\ref{prop:maxmon-resolvent-Yosida}, item (vi).
\end{rem}
}
\subsection{Relevance of Theorem~\ref{thm:main}}
In what follows we will show how the results of \cite{KunzMont97} as well as the results on sweeping processes can be recovered from Theorem~\ref{thm:main}.

\subsubsection{Recovering the results of \cite{KunzMont97}}\label{sec:oscac}
As recalled in Sect.~\ref{sec:connections}, the results in~\cite{KunzMont97} imposed continuity with respect to the pseudo-distance introduced in Definition~\ref{def:dis}. We claim that if the mapping $t\mapsto F(t,\cdot):\R^n\rightrightarrows\R^n$ is absolutely continuous on $[0,T]$ then the set-valued mapping $t\mapsto \graph F(t,\cdot)$ is {outer semicontinuous from right} on $[0,T]$. To see this, let $(t_\ell,x_\ell,y_\ell)_{\ell\in N}\subseteq[0,T]\times\R^n\times\R^n$ be a sequence such that $y_\ell\in F(t_\ell,x_\ell)$, {$t_{\ell} \ge t$},  and $\lim_{\ell\uparrow\infty}(t_\ell,x_\ell,y_\ell)=(t,x,y)$ for some $t\in[0,T]$, $x,y\in\R^n$. What needs to be proven is that $y\in F(t,x)$. To see this, let $(\eta,\zeta)\in\R^{n+n}$ be such that $\zeta\in F(t,\eta)$. From absolute continuity of $t\mapsto F(t,\cdot)$, we have
$$
\inn{\zeta-y_\ell}{x_\ell-\eta}\leq \big(r(t)-r(t_\ell)\big)\big(1+\abs{\zeta}+\abs{y_\ell}\big).
$$
By letting $\ell$ tend to infinity, we obtain
$$
\inn{\zeta-y}{x-\eta}\leq 0
$$
since $r$ is continuous. This means that $y\in F(t,x)$ as $F(t,\cdot)$ is maximal monotone. Another hypothesis required by the results of \cite{KunzMont97} is a linear growth condition that coincides with \ref{ass:least}.

{
\subsubsection{Inclusions with subdifferentials of convex functions}\label{sec:convexSubdiff}
Another closely related class of dynamical systems appearing in the literature is given by
\begin{equation}\label{eq:subdiffConvex}
\dot x \in -\partial g(t,x), \quad x(0) = x_0 \in \dom \partial g(0,\cdot),
\end{equation}
where $g(t,\cdot):\R^n \to \R \cup \{\infty\}$ is a convex, lower semicontinuous, and proper for each $t \ge 0$, and $\partial g(t,\cdot)$ denotes its subdifferential.
\begin{cor}\label{cor:convexSubdiff}
Consider the differential inclusion \eqref{eq:subdiffConvex}, and assume that the following hold:
\begin{itemize}
\item There exists a nondecreasing function $\varphi\in AC([0,T],\R)$ such that 
\[
 \sup_{z \in \dom \partial g(s,\cdot)}\dist \big(z , \dom \partial g (t,\cdot)\big)  \leq \varphi(t)- \varphi(s), \quad \forall s,t \text{ with } 0\leq s\leq t \leq T.
 \]
\item There exist a nondecreasing function $\theta:[0,T]\rightarrow\R$, and a scalar $\Lambda > 0$,  such that for every $0 < \lambda < \Lambda $,
\[
\abs{\nabla g_\lambda(s+\lambda,x)} \leq \theta(s+\lambda)-\theta(s)+\abs{\nabla g_\lambda(s,x)}
\]
for all $s \text{ with } 0\leq s \leq T - \lambda$, and $x\in \dom \partial g(s,\cdot)$, where $g_{\lambda}(t,\cdot)$ denotes the continuously differentiable Yosida approximation of $g(t,\cdot)$, defined as
\[
g_\lambda(t,x) = \inf_{y \in \R^n} \left \{ g(t,y) + \frac{1}{2\lambda} \|y -x \|^2 \right\}.
\]
\end{itemize}
Then, the differential inclusion \eqref{eq:subdiffConvex} admits a unique solution.
\end{cor}

The corollary, stated above, is a direct consequence of Theorem~\ref{thm:main}. Under the assumptions on $g(t,\cdot)$, \ref{ass:basicF} holds. Moreover, Yosida approximation of the operator $\partial g$ can be expressed directly in terms of the Yosida approximation of the function $g$ itself. We also have the following result which ensures that, in such cases, \ref{ass:dist} implies \ref{ass:new}.

\begin{lem}\label{lem:A3ImpliesA2}
If, for each $t \ge 0$, $F(t,\cdot) = \partial g(t,\cdot)$ for some $g(t,\cdot)$ convex, lower semicontinuous and proper. Then, \ref{ass:dist} implies \ref{ass:new}.

\BP
Let $(t_\ell,x_\ell,y_\ell)_{\ell\in \N}\subseteq[0,T]\times\R^n\times\R^n$ be a sequence such that $y_\ell\in \partial g(t_\ell,x_\ell)$ and moreover, $\lim_{\ell\uparrow\infty}(t_\ell,x_\ell,y_\ell)=(t,x,y)$, $t_\ell \geq t$, for some $t\in[0,T]$, $x,y\in\R^n$. We now show that $y\in \partial g(t,x)$. For that, let $(\eta,\zeta)\in\R^{n+n}$ be such that $\zeta\in \partial g(t,\eta)$, and let $x_{\ell,\eta}^\star := \argmin_{z \in \dom \partial g(t_\ell,\cdot)} \vert z - \eta\vert$. Note that, $\dom \partial g(t,\cdot)$ is closed, convex and nonempty for each $t \ge 0$, and hence $x_{\ell,\eta}^\star$ is unique, and because of \ref{ass:dist}, $x_{\ell,\eta}^\star \to \eta$ as $\ell \to \infty$. We, therefore, have
\begin{align*}
\inn{y_\ell}{\eta - x_\ell} & =  \inn{y_\ell}{\eta - x^\star_{\ell,\eta}} + \inn{y_\ell}{x^\star_{\ell,\eta} - x_\ell} \leq g(t_\ell, x_{\ell,\eta}^\star) - g(t_\ell, x_{\ell}) + \inn{y_\ell}{\eta - x^\star_{\ell,\eta}} \\
& \leq g(t_\ell, x_{\ell,\eta}^\star) - g(t_\ell, x_{\ell}) + \vert y_\ell \vert \ \dist (\eta, \dom \partial g(t_\ell,\cdot)) \\
& \leq g(t_\ell, x_{\ell,\eta}^\star) - g(t_\ell, x_{\ell}) + \vert y_\ell \vert \ \dist ( \dom \partial g(t,\cdot), \dom \partial g(t_\ell,\cdot) ) \\
& \leq g(t_\ell, x_{\ell,\eta}^\star) - g(t_\ell, x_{\ell}) + \vert y_\ell \vert \ (\varphi(t_\ell) - \varphi(t)) .
\end{align*}
By letting $\ell$ tend to infinity, and invoking \ref{ass:dist}, we obtain
\[
\inn{y}{\eta - x}  \le g(t, \eta) - g(t, x) 
\] 
which implies $y \in \partial g(t,x)$, and hence assumption \ref{ass:new} holds.
\EP
\end{lem}

\begin{example}
As an illustration of Corollary~\ref{cor:convexSubdiff}, consider the function $g:\R_{+} \times \R \to \R_{+}$, defined by, $g(t,x) = t \abs{ x }$. The first item holds since the domain of $\partial g$ does not change with time. To check the second condition, Yosida approximation of the function $g(t,\cdot)$, for each $t \ge 0$, is described by:
\[
g_{\lambda}(t,x) = \begin{cases} \frac{1}{2\lambda} \abs{x}^2 & \text{ if } \abs{x} \le t \lambda \\ t \abs{x} - \frac{\lambda t^2}{2} & \text{ if } \abs{x} > t \lambda \end{cases}.
\]
As a result, Yosida approximation of the operator $\partial g(t, \cdot)$, for each $t \ge 0$, is given by:
\[
\partial g_{\lambda}(t,x) = \begin{cases} \frac{1}{\lambda} x & \text{ if } \abs{x} \le t \lambda \\ t \ {\rm sign}(x) & \text{ if } \abs{x} > t \lambda \end{cases}
\]
and hence it follows that for each $t \ge s \ge 0$, and each $x \in \R$, we have
\[
\vert \partial g_{\lambda}(t,x) \vert - \vert \partial g_{\lambda}(s,x) \vert \le t -s
\]
and hence \ref{ass:static2} holds with linear function $\theta(s) = s$.
\end{example}
}

\subsubsection{Recovering the special case of sweeping processes}\label{sec:sweep}
Sweeping process is a special case of \eqref{eq:incGen}, where $F(t,x) = \cN_{S(t)}(x)$, with $S:[0,T] \rightrightarrows \R^n$ being a closed convex-valued mapping. The normal cone operator is by definition maximal monotone for each $t \in [0,T]$. Moreover, in this particular case, assumption~\ref{ass:least} is trivially satisfied since $\cN_{S(t)}(x)$ is a cone and thus $0 \in \cN_{S(t)}(x)$, for each $t \in [0,T]$ and $x \in \dom \, \cN_{S(t)}(\cdot) = S(t)$. {Assumption \ref{ass:dist} can be written directly in terms of the set-valued mapping $S$ as follows: There exists a nondecreasing function $\varphi \in AC([0,T],\R)$ such that
\[
\dist(S(s),S(t)) :=\sup_{z \in S(s)} \dist (z,S(t)) \leq \varphi(t)- \varphi(s), \quad \forall s,t \text{ with } 0\leq s\leq t \leq T.
\]
The aforementioned condition is equivalent to saying that the retraction of the set-valued mapping $S(\cdot)$ is absolutely continuous, and it was first proposed in \cite{More77} for existence of solutions to first order sweeping processes.
One can invoke Lemma~\ref{lem:A3ImpliesA2} to check that assumption \ref{ass:dist} implies \ref{ass:new}. Indeed, by taking $g(t,x) = \begin{cases}0 & \text{if } x \in S(t), \\ +\infty & \text{if } x \not\in S(t) \end{cases}$, and observing that $\cN_{S(t)}(\cdot) = \partial g{(t,\cdot)}$.

\begin{rem}
{In Section~\ref{sec:convexSubdiff}, and Section~\ref{sec:sweep}, we considered a specific class of systems where the maximal monotone operator is defined by the subdifferential of the convex function. In the literature, we also find other results for differential inclusions with subdifferentials of time-dependent convex functions, see for example \cite[Page 197]{HuPapa00}, \cite{Kubo88}, \cite{Papa90}, and \cite{Tols17}. In some of these works, conditions are stated in terms of the function $g(\cdot,\cdot)$ itself, whereas we work with assumptions that result in the bounds on the corresponding Yosida approximations of the subdifferentials. Some results from \cite{Kubo88} can be used to write the conditions in terms of the function $g(\cdot,\cdot)$ itself, but in general, it is not clear how our conditions on the subdifferential operator can be translated to the conditions on the function $g(\cdot,\cdot)$.
}
\end{rem}
}

\section{Proof of Theorem~\ref{thm:main}}\label{sec:proof}
We are basically concerned with the existence of the solution in this proof, as the uniqueness readily follows from assumption \ref{ass:basicF}.
The proof of existence is based on constructing a sequence of approximate solutions and showing that this sequence converges to a function which satisfies the differential inclusion \eqref{eq:incGen}. This is formally done in following main steps:
\begin{itemize}
\item Discretizing \eqref{eq:incGen}
\item Obtaining bounds on discrete values
\item Construction of a sequence of approximate solutions
\item Studying the limit of the sequence
\end{itemize}
Each of these steps is carried out as a subsection in the sequel, and it is shown that the limit we thus obtain is indeed a solution to \eqref{eq:incGen}.

\subsection{Discretization of \eqref{eq:incGen}}
We first begin with discretizing \eqref{eq:incGen}. Let
$$
\Delta=\set{t_0,t_1,\ldots,t_K}{0 = t_0 < t_1 < \cdots < t_k < t_{k+1} < \cdots < t_K=T}
$$
be a partition of the interval $[0,T]$. Define
$$
0 < h_k:= t_{k}-t_{k-1}
$$
for $k\in\nset{K}$. Note that $\sum_{k=1}^{K} h_k=T$. We define the {\em size\/} of the partition $\Delta$ by $K(\Delta)$ and the {\em granularity\/} by $\abs{\Delta}=\max_{k\in\nset{K}} h_k$. For simplicity, we write $K=K(\Delta)$ when $\Delta$ is clear from the context.

Next, consider the discretization of \eqref{eq:incGen} based on the partition $\Delta$ given by
\beq\label{eq:iteration}
\frac{x_{k+1}-x_k}{h_{k+1}}\in -F(t_{k+1},x_{k+1})
\eeq
for $k\in\pset{0,1,\ldots,K-1}$. Alternatively, we have
\beq\label{e:discrete-yosida}
x_{k+1} = \big( I +h_{k+1}  F (t_{k+1},\cdot)\big)^{-1}(x_{k}).
\eeq
This resolvent  based alternative form, together with assumption \ref{ass:basicF} and Proposition~\ref{prop:maxmon-resolvent-Yosida}, guarantees that the discretization \eqref{eq:iteration} is well-defined in the sense that there exist $x_0,x_1,\ldots,x_K$ satisfying \eqref{eq:iteration} (and hence \eqref{e:discrete-yosida}). If $\dom(F(t,\cdot))$ is closed for each $t \geq 0$, then it follows from \eqref{e:discrete-yosida} that $x_{k+1}$ is a projection of $x_k$ on the set $\dom(F(t_{k+1},\cdot))$; Put simply, in case the domain of the multivalued function $F(t,\cdot)$ is closed, the sequence of points $x_k$ is obtained by applying the proximal operator associated with the mapping $h_k F(t_k,\cdot)$, as defined in \cite{Rock76}, \cite{PariBoyd13}. We will use a certain interpolation between the points $x_k$ to get an approximate solution for the differential inclusion~\eqref{eq:incGen}.

%

\subsection{Bounds on $x_k$ values}
{We aim at establishing bounds on $x_k$ that are independent of the underlying partition $\Delta$. The bounds obtained under the condition \ref{ass:least} and \ref{ass:static} are different and appear in Lemma~\ref{lem:bndSeq} and Lemma~\ref{lem:bndSeqYB}, respectively. We prefer using the same notation for the bounds of our interest in these lemmas, and when we use these bounds later on, the reference to these bounds will be clear from the context. 
}

\begin{lem}\label{lem:bndSeq}
Under assumption \ref{ass:least}, there exists a sequence of partitions $\{\Delta_\ell\}$ with $\vert \Delta_\ell \vert \to 0$, such that for each of these partitions, we have
\begin{align}
\abs{x_k}&\leq \beta \label{i:uni bound}\\
\abs{x_{k} - x_{k-1}}  & \leq \psi(t_k) - \psi(t_{k-1})\label{e:bound with psi}
\end{align}
for each $k\in\nset{K}$; the constant $\beta \ge 0$, and the function $\psi:[0,T]\rightarrow\R_+$ are defined as
\begin{gather}
\beta=\alpha+\varphi(T)-\varphi(0)+(1+\alpha)\int_0^T\!\!\sigma(s)\,ds\label{e:only beta}\\
\psi(t):=t+2\varphi(t)+ (1+\gamma)\int_0^t \sigma(s)\,ds\quad\forall\,t\in[0,T], \label{eq:defPsi}
\end{gather}
with $\varphi$ satisfy \ref{ass:dist}, $\sigma$ satisfying \ref{ass:least}, and
\begin{gather}
\alpha=\abs{x_0}+\varphi(T)-\varphi(0)\label{e:alpha and r alpha} \\
\gamma=\beta+\varphi(T)-\varphi(0). \label{e:gamma and r gamma}
\end{gather}
\end{lem}

\BP 
To obtain the bounds given in \eqref{i:uni bound} and \eqref{e:bound with psi}, we start by analyzing the sequence \eqref{e:discrete-yosida} for a fixed partition and introduce some simplified notation for the corresponding operators:
\begin{equation}\label{eq:notaDisc}
F_k:=F (t_{k},\cdot),\quad
J_k:=\big(I +h_{k} F_k \big)\inv,
\qand
Y_k:=\frac{1}{h_k}(I-J_k).
\end{equation}
It then follows from \eqref{e:discrete-yosida} that
\begin{gather}
x_{k+1}\in \dom F_{k+1}\label{e:x k+1 props.1}\\
x_{k+1} = J_{k+1}(x_{k})\label{e:x k+1 props.2}
\end{gather}
for all $k\in\pset{0,1,\ldots,K-1}$, where we recall that $K$ is the size of the chosen partition.

To establish \eqref{i:uni bound}, we first introduce auxiliary points $\bar x_k$ given by $\bar x_0=x_0$ and
\[
\bar x_{k+1} := \proj\big(\bar x_{k}, \cl(\dom{ F_{k+1}})\big)
\]
for $k\in\pset{0,1,\ldots,K-1}$.
Clearly, we have
\beq\label{e:overline x in dom}
\bar x_k\in\cl(\dom F_{k})
\eeq
for all $k\in\pset{0,1,\ldots,K}$. Then, it follows from assumption \ref{ass:dist} and \eqref{e:hausdorff conseq} that
\beq\label{eq:boundxbarder}
\abs{\bar x_{k} - \bar x_{k-1}} \leq\varphi(t_{k}) - \varphi(t_{k-1})
\eeq
for all $k\in\pset{1,2,\ldots,K}$. We thus obtain, for each $k\in\nset{K}$,
\bse
\begin{align}
	|\bar x_{k}| &\leq \abs{\bar x_{k-1}}+\abs{\bar x_{k} - \bar x_{k-1}} \leq|\bar x_{k-1}| + \varphi(t_{k}) - \varphi(t_{k-1})\label{e:bar x.1} \\
	& \leq\vert \bar x_0 \vert + \varphi(t_{k}) - \varphi (t_0)\label{e:bar x.2} \\
	& \leq \abs{\bar x_0}+\varphi(T) - \varphi (0)\label{e:bar x.3}
\end{align}
\ese
where \eqref{e:bar x.2} follows from the repeated application of \eqref{e:bar x.1}, and \eqref{e:bar x.3} uses the fact that $\varphi$ is nondecreasing. An immediate consequence of \eqref{e:bar x.3} is that 
\beq\label{eq:bound}
	\abs{ \bar x_k} \leq \alpha
\eeq
for each $k\in\pset{0,1,\ldots,K}$ where $\alpha$ satisfies \eqref{e:alpha and r alpha}.

The bounds given in \eqref{i:uni bound} and \eqref{e:bound with psi} are now obtained from assumption~\ref{ass:least}.
From \eqref{e:overline x in dom}, we have that, for every $\varepsilon>0$ and each $k \in \{1,2,\dots, K\}$, there exists a point $\bar x_k^\varepsilon$ satisfying
\begin{gather}
\bar x_k^\varepsilon\in\dom F_k\label{e:bar x eps in dom}\\
\abs{\bar x_k^\varepsilon-\bar x_k}\leq \varepsilon , \label{e:from bar x to bar x eps}
\end{gather}
and therefore
\beq\label{e:diff bar x eps and bar x}
\abs{\bar x_k^\varepsilon}\leq \alpha+\varepsilon
\eeq
for all $k\in\pset{0,1,\ldots,K}$ in view of \eqref{eq:bound} and the triangle inequality. Next, we introduce the sequence of points $\bar y_k^\varepsilon$ based on $\bar x_k^\varepsilon$ by
\begin{align}\label{eq:resolventproj}
	\bar y_{k}^\varepsilon := J_k(\bar x_{k}^\varepsilon)
\end{align}
for all $k\in\pset{0,1,\ldots,K}$.
Note that
\beq
\dfrac{\bar x_{k}^\varepsilon-\bar y_{k}^\varepsilon}{h_{k}} = Y_{k}(\bar x_{k}^\varepsilon).
\eeq
Now, it follows from \eqref{e:bar x eps in dom} and Proposition~\ref{prop:maxmon-resolvent-Yosida} that
\beq\label{eq:discDiffBnd}
	\left | \dfrac{\bar x_{k}^\varepsilon - \bar y_{k}^\varepsilon}{h_{k}} \right | \leq | F ^0_k(\bar x_{k}^\varepsilon)|.
\eeq
To obtain a bound on the right-hand side of \eqref{eq:discDiffBnd}, we employ assumption \ref{ass:least}. Take any partition $\Delta_\ell$, and assume momentarily that $\sigma$ is a constant function.

By using assumption \ref{ass:least} and \eqref{e:diff bar x eps and bar x}, we get
\beq\label{e:bound uni in eps}
\abs{ \bar x_{k}^\varepsilon - \bar y_{k}^\varepsilon} \leq  h_k \sigma (1+\alpha+\varepsilon)
\eeq
for all $k\in\pset{0,1,\ldots,K}$ and $\eps > 0$. In view of \eqref{e:x k+1 props.2} and \eqref{eq:resolventproj}, Proposition~\ref{prop:maxmon-resolvent-Yosida} implies that
\beq\label{e:from x k to x k-1}
	|x_{k}-\bar y_{k}^\varepsilon| \leq |x_{k-1}-\bar x_{k}^\varepsilon|.
\eeq
Hence, using \eqref{e:bound uni in eps} and \eqref{e:from x k to x k-1}, we obtain
\beq\label{eq:diffproj}
	|x_{k}-\bar x_{k}^\varepsilon| \leq |x_{k}-\bar y_{k}^\varepsilon| + |\bar y_{k}^\varepsilon-\bar x_{k}^\varepsilon| \leq
	|x_{k-1}-\bar x_{k}^\varepsilon|+h_k \sigma (1+\alpha+\varepsilon)
\eeq
for all $k\in\nset{K}$ and $\eps > 0$. Letting $\eps$ tend to zero in \eqref{eq:diffproj} results in
\begin{align}\label{eq:diffproj2}
	\abs{x_{k}-\bar x_{k}} \leq \abs{x_{k-1} - \bar x_{k}} + h_k \sigma(1+\alpha)
\end{align}
for all $k\in\nset{K}$. Therefore, we have
\begin{align}
	|x_{k} - \bar x_{k}| & \leq \abs{x_{k-1} - \bar x_{k}} + h_k \sigma (1+\alpha)\\
	& \leq|x_{k-1} - \bar x_{k-1}| +\abs{\bar x_{k-1}- \bar x_{k}} + h_k \sigma (1+\alpha) \\
	&\stackrel{\eqref{eq:boundxbarder}}{\leq}  |x_{k-1} - \bar x_{k-1}| + \varphi(t_k) - \varphi(t_{k-1})+ h_k \sigma(1+\alpha)\\
	& \leq|x_0 - \bar x_0|+ \varphi(t_k) - \varphi(0) + (\sum_{\ell=1}^{k}h_k) \sigma (1+\alpha) \\
	&\leq \varphi(T) - \varphi(0)+T \sigma (1+\alpha)\label{eq:diffxbarx}.
\end{align}
for all $k\in\nset{K}$. Since the constant $\sigma$ is finite, we can conclude from \eqref{eq:bound} and \eqref{eq:diffxbarx} that 
\begin{align}\label{eq:boundx}
	\abs{x_{k}} \leq \beta
\end{align}
for all $k\in\nset{K}$ where $\beta$ satisfies \eqref{e:only beta}. This establishes \eqref{i:uni bound}.

Next, we proceed to establish \eqref{e:bound with psi} using the bound in assumption \ref{ass:least} and \eqref{i:uni bound}. 
To this end, we continue using the notation introduced in \eqref{eq:notaDisc}, and introduce a sequence of auxiliary points $\xi_{k}$
\begin{align}
	\xi_{k}:=\proj{\big(x_{k-1},\cl(\dom  F_k)\big)}
\end{align}
for all $k\in\nset{K}$. Clearly, we have
\beq\label{e:xi in dom}
\xi_{k}\in\cl(\dom F_{k})
\eeq
for all $k\in\nset{K}$. Therefore, for all $\varepsilon>0$ there exist points $\xi_{k}^\varepsilon$ satisfying
\bse
\begin{gather}
\xi_{k}^\varepsilon\in\dom F_k\label{e:xi eps in dom}\\
\abs{\xi_{k}^\varepsilon-\xi_k}\leq \varepsilon\label{e:from xi to xi eps}.
\end{gather}
\ese
It follows from \eqref{e:hausdorff conseq} and assumption~\ref{ass:dist} that
\begin{align}\label{eq:xxbarbound}
	|x_{k-1}-\xi_{k}| \leq \varphi(t_k) - \varphi(t_{k-1})
\end{align}
for all $k\in\nset{K}$. In view of \eqref{eq:boundx} and the fact that $\varphi$ is nondecreasing, this means that
$$
\abs{\xi_{k}}\leq \beta+\varphi(T)-\varphi(0)
$$
for all $k\in\nset{K}$. As such, we have
\beq\label{e:xi k eps}
\abs{\xi^\eps_k}\leq \gamma+\eps
\eeq
for all $k\in\nset{K}$ where $\gamma$ satisfies \eqref{e:gamma and r gamma}. Now, define
\begin{align}
	\zeta_{k}^\varepsilon :=J_k(\xi_{k}^\varepsilon)
\end{align}
for all $k\in\nset{K}$. Note that
\begin{align}
\vert x_{k} - x_{k-1} \vert & \leq\vert x_k - \zeta_k^\varepsilon \vert + \vert \zeta_k^\varepsilon - \xi_k^\varepsilon \vert + \vert \xi_k^\varepsilon - x_{k-1} \vert \notag\\
& \leq2 \vert \xi_k^\varepsilon - x_{k-1} \vert + \vert \zeta_k^\varepsilon - \xi_k^\varepsilon \vert \label{eq:bndDiff}
\end{align}
where we used the fact that $\vert x_k- \zeta_k^\varepsilon \vert \leq\vert x_{k-1}-\xi_k^\varepsilon\vert$ due to the resolvent  being nonexpansive. 
Since $\xi_k^\varepsilon \in \dom F_k$ due to \eqref{e:xi eps in dom}, we can invoke the bound on the least norm element of $ F (t_k, \xi_k^\varepsilon)$ in assumption \ref{ass:least} to obtain
\[
\frac{\xi_k^\varepsilon - \zeta_k^\varepsilon}{h_k} = Y_{k}(\xi_k^\varepsilon) \Longrightarrow \left \vert \xi_k^\varepsilon - \zeta_k^\varepsilon \right \vert \leq h_k \,  \vert  F ^0(t_k,\xi_k^\varepsilon) \vert \overset{\eqref{e:xi k eps}}{\leq} h_k \, \sigma (1+\gamma+\eps)
\]
for all $k\in\pset{1,\ldots,K}$ and $\eps > 0$. Together with \eqref{e:from xi to xi eps} and \eqref{eq:xxbarbound},
\eqref{eq:bndDiff} leads to
\[
\vert x_{k} - x_{k-1} \vert \leq 2 \big(\varphi(t_k) - \varphi(t_{k-1})\big) + h_k \sigma(1+\gamma)
\]
by taking the limit as $\varepsilon$ tends to zero. This establishes the bound given in \eqref{e:bound with psi} for constant $\sigma$ and arbitrary partition. 

To obtain the bounds \eqref{i:uni bound} and \eqref{e:bound with psi} for $\sigma$ locally integrable, and not necessarily constant, one can invoke \cite[Lemma 3.3.1]{LaksLeel81}. To use this result, one must choose the partitions $\Delta_\ell$ carefully, and work with a piecewise constant approximation (in $L_1$-norm) of the function $\sigma$. 
\EP


{
\begin{lem}\label{lem:bndSeqYB}
Under the assumption \ref{ass:static}, for any partition $\Delta$, with $\abs{\Delta} < \Lambda$ we have
\begin{align}
\abs{x_k}&\leq \beta \label{i:uni bound with YB} \\
\abs{x_{k} - x_{k-1}}  & \leq \psi(t_k) - \psi(t_{k-1}) \label{e:bound with psi with YB}
\end{align}
for each $k\in\nset{K}$; the constant $\beta \ge 0$, and the function $\psi : [0,T]\rightarrow\R_+$ are defined as
\begin{gather}
\beta :=\abs{x_0} + \alpha T \label{eq:betaYB} \\
\psi(t) :=\alpha t, \label{eq:defPsiYB}
\end{gather}
with $\alpha:= 2(\theta(T) -\Theta(0)) + F^0(0,x_0)$.
\end{lem}

\BP
Let $x_0 \in \dom F(0,\cdot)$. 
%
Since $x_{k+1} = J_{k+1}(x_k) \in \dom F_{k+1}$, we observe that
\begin{equation}\label{eq:tempDiff}
\abs{\frac{x_{k+1} - x_k}{h_{k+1}}} = \abs{Y_{k+1}(x_k)} \le \theta(t_{k+1}) - \theta(t_k) + \abs{F^0_{k}(x_k)}.
\end{equation}
On the other hand,
\[
Y_{k+1}(x_k) \in F_{k+1}(J_{k+1}(x_k)) 
\]
and hence,
\[
\abs{Y_{k+1}(x_k)} \ge \abs{F_{k+1}^0(x_{k+1})} . 
\]
We therefore have
\begin{align*}
\abs{F^0_{k+1}(x_{k+1})} \le \theta(t_{k+1}) - \theta(t_k) + \abs{F^0_{k}(x_k)}
\end{align*}
and, in particular, for each $k \in \N$,
\begin{align*}
\abs{F^0_{k}(x_{k})} \le \theta(t_{k+1}) - \theta(0) + \abs{F^0_{0}(x_0)}.
\end{align*}
Plugging this bound in \eqref{eq:tempDiff}, and using the fact that $\theta$ is nondecreasing, we obtain
\[
\abs{\frac{x_{k+1} - x_k}{h_{k+1}}} \le 2(\theta(T) - \theta(0)) + \abs{F^0_{0}(x_0)}.
\]
Letting $\alpha:= 2(\theta(T) - \theta(0)) + F_0^0(x_0)$, we get
\[
\abs{x_{k+1} - x_k} \le \alpha (t_{k+1} - t_k),
\]
so that \eqref{e:bound with psi with YB} holds with $\psi(s) = \alpha \, s$. By repeated application of the trying inequality, $\abs{x_k} \le \abs{x_{k-1}} + \abs{x_{k} - x_{k-1}}$, we get
\[
\abs{x_k} \le \abs{x_0} + \alpha T,
\]
so that \eqref{i:uni bound with YB} holds with $\beta :=\abs{x_0} + \alpha T$.
\EP

In subsequent analysis, we refer to $\beta$ and $\psi$, either from Lemma~\ref{lem:bndSeq} or Lemma~\ref{lem:bndSeqYB}, depending upon whether we assume \ref{ass:least} or \ref{ass:static}.
}
\subsection{Construction of a sequence of approximate solutions}
Based on the $x_k$ values, we construct a sequence of absolutely continuous (in time) functions which approximate the actual solution of the system. To this end, note that the function $\psi$ defined above is strictly increasing and absolutely continuous. Now, define the piecewise continuous function $x_\Delta$ as
\beq\label{e:x delta}
x_\Delta(t) :=
\dfrac{\psi(t_{k+1})-\psi(t)}{\psi(t_{k+1})-\psi(t_k)}x_k+\dfrac{\psi(t)-\psi(t_k)}{\psi(t_{k+1})-\psi(t_k)}x_{k+1}
\eeq
where $t \in [t_{k},t_{k+1}]$ and $k\in\pset{0,1,\ldots,K-1}$. By definition, $x_\Delta$ is a continuous function and
\beq\label{e:x delta at tk}
x_\Delta(t_k)=x_k
\eeq
for all $k\in\pset{0,1,\ldots,K}$. We will show that
\[
x(t) := \lim_{\vert \Delta \vert \to 0} x_{\Delta} (t)
\]
is the desired solution to the inclusion \eqref{eq:incGen}. An important intermediate step in studying the convergence of the sequence $x_{\Delta}$ is to obtain the following uniform bound.

\begin{lem}\label{lem:bndEquiCont}
Let $\ul\tau$ and $\ol \tau$ be such that $0\leq \ul\tau < \ol \tau\leq T$. For any partition $\Delta$, it holds that
\beq\label{e:towards equi cont}
\abs{x_\Delta(\ol \tau)-x_\Delta(\ul\tau)}\leq \psi(\ol\tau) - \psi(\ul\tau).
\eeq
\end{lem}
\BP From the definition of $x_{\Delta}$ \eqref{e:x delta} for a fixed partition $\Delta$,
there exist integers $q$ and $r$ with $q+1\leq r$ such that $t_q\leq \ul\tau < t_{q+1}$ and $t_{r-1}<\ol\tau\leq t_{r}$. If $q+1=r$, then we have
\begin{alignat}{3}
\abs{x_\Delta(\ol\tau)-x_\Delta(\ul\tau)}
&\leq  \abs{\frac{\psi(\ol\tau)-\psi(\ul\tau)}{\psi(t_{q+1})-\psi(t_q)}(x_{q+1}-x_q)}
&\qquad&\text{from \eqref{e:x delta}} \notag\\
&\leq \frac{\psi(\ol\tau)-\psi(\ul\tau)}{\psi(t_{q+1})-\psi(t_q)}\big(\psi(t_{q+1}) - \psi(t_{q})\big)
&\qquad&\text{from \eqref{e:x delta at tk}}\notag\\
&\leq \psi(\ol\tau)-\psi(\ul\tau).
&\qquad&\label{e:towards equi cont1}
\end{alignat}
In a similar fashion, if $q+1<r$ then we have
\begin{alignat}{3}
\abs{&x_\Delta(\ol\tau)-x_\Delta(\ul\tau)}  \leq \abs{x_\Delta(t_{r-1})-x_\Delta(\ol\tau)} \notag 
+ \sum_{q+1\leq i\leq r-2}\abs{x_\Delta(t_{i+1})-x_\Delta(t_i)}+\abs{x_\Delta(\ul\tau)-x_\Delta(t_{q+1})} \notag\\
&=\frac{\psi(\ol\tau)-\psi(t_{r-1})}{\psi(t_{r})-\psi(t_{r-1})}\abs{x_{r}-x_{r-1}} + \sum_{q+1\leq i\leq r-2}\abs{x_\Delta(t_{i+1})-x_\Delta(t_i)} \notag + \frac{\psi(t_{q+1})-\psi(\ul\tau)}{\psi(t_{q+1})-\psi(t_q)}\abs{x_{q+1}-x_{q}}\notag\\
& \leq \psi(\ol\tau)-\psi(t_{r-1})  + \sum_{q+1\leq i\leq r-2}\big(\psi(t_{i+1})-\psi(t_i)\big)+\psi(t_{q+1})-\psi(\ul\tau) =  \psi(\ol\tau)-\psi(\ul\tau).
\end{alignat}
Hence, \eqref{e:towards equi cont} is established.\EP

\subsection{Limit of the sequence}
The bounds established in the previous section allow us to study the limiting behaviour of the sequence $(x_{\Delta_\ell})_{\ell \in \N}$.

\begin{lem}\label{lem:equiCont}
Consider a sequence of partitions $(\Delta_\ell)_{\ell\in\N}$ with $\abs{\Delta_\ell}\rightarrow 0$ as $\ell$ tends to infinity. The sequence $(x_{\Delta_\ell})_{\ell\in\N}$ is equicontinuous.
\end{lem}

\BP
Note first that $\psi$ introduced in \eqref{eq:defPsi} (or, \eqref{eq:defPsiYB}) is uniformly continuous on the compact interval $[0,T]$ as it is absolutely continuous on the same interval.  Therefore, for any $\varepsilon>0$ there exists a positive number $\delta>0$ such that
$$
\abs{\psi(\ol\tau)-\psi(\ul\tau)}<\varepsilon
$$
for all $\ul\tau,\ol\tau\in[0,T]$ such that $\abs{\ol\tau-\ul\tau}<\delta$. In view of \eqref{e:towards equi cont}, we have
$$
\abs{x_{\Delta_\ell}(\ol\tau)-x_{\Delta_\ell}(\ul\tau)}<\varepsilon
$$
for all $\ell\in\N$ and $\ul\tau,\ol\tau\in[0,T]$ such that $\abs{\ol\tau-\ul\tau}<\delta$. Consequently, the sequence $(x_{\Delta_\ell})_{\ell\in\N}$ is equicontinuous.
\EP

Let $(\Delta_\ell)_{\ell\in\N}$ be a sequence of partitions with $\abs{\Delta_\ell}\rightarrow 0$ as $\ell$ tends to infinity. Since the sequence $(x_{\Delta_\ell})_{\ell\in\N}$ is also uniformly bounded in view of Lemma~\ref{lem:bndSeq} (as well as Lemma~\ref{lem:bndSeqYB}), Theorem~\ref{thm:ascArz} (Arzel\'{a}-Ascoli theorem) implies that it converges uniformly to a continuous function $x$ on a subsequence, say $N\in\calN_\infty^\#$. We claim that $x$ is absolutely continuous. To see this, let $\ul\tau,\ol\tau\in[0,T]$ with $\ul\tau\leq \ol\tau$ and
note that
\begin{align}
\abs{x(\ol\tau)-x(\ul\tau)}&\leq\abs{x(\ol\tau)-x_{\Delta_\ell}(\ol\tau)}+\abs{x_{\Delta_\ell}(\ol\tau)-x_{\Delta_\ell}(\ul\tau)}+\abs{x_{\Delta_\ell}(\ul\tau)-x(\ul\tau)}\notag\\
&\leq\abs{x(\ol\tau)-x_{\Delta_\ell}(\ol\tau)}+\psi(\ol\tau)-\psi(\ul\tau)+\abs{x_{\Delta_\ell}(\ul\tau)-x(\ul\tau)}\notag\\
&\leq \psi(\ol\tau)-\psi(\ul\tau)\label{e:bound x with psi}
\end{align}
where the first inequality follows from the triangle inequality, the second from \eqref{e:towards equi cont}, and the third by taking the limit on the convergent subsequence $N$. Thus, absolute continuity of $x$ follows from absolute continuity of the function $\psi$.

Now, we want to show that $x$ is a solution of \eqref{eq:incGen}, that is
\beq\label{e:x is a solution or not}
x(t) \in \dom F(t,\cdot) \qand \dot{x}(t)\in -F\big(t,x(t)\big)
\eeq
for almost all $t\in [0,T]$.

Let $\Gamma\subseteq[0,T]$ be defined by $\Gamma=\set{t\in(0,T)}{\psi\text{ and }x\text{ are both }\text{differentiable at } t\text{ and }t\not\in\cup_{\ell\in N}\Delta_\ell}$. Since $\psi$ and $x$ are both absolutely continuous and $\cup_{\ell\in N}\Delta_\ell$ is countable, it is enough to show \eqref{e:x is a solution or not} for almost all $t\in\Gamma$.

For a partition $\Delta$, define 
\beq
y_{\Delta}(t)=\frac{x_{k+1}-x_{k}}{\psi(t_{k+1})-\psi(t_k)}
\eeq
for $t\in(t_k,t_{k+1})$ and $y_{\Delta}(t_k)=0$ for $t_k\in\Delta$. 

From \eqref{e:x delta}, we see that
\beq\label{e:x delta ell deriv}
\dot{x}_{\Delta_\ell}(t)=\dot{\psi}(t)\frac{x_{k+1}-x_{k}}{\psi(t_{k+1})-\psi(t_k)}=\dot{\psi}(t)y_{\Delta_\ell}(t)
\eeq
for all $t\in\Gamma$. 

In view of \eqref{e:x delta} and Lemma~\ref{lem:bndEquiCont}, we see that $\abs{y_{\Delta_\ell}}_{\Linf}\leq 1$ for all $\ell$. Therefore, the sequence $(y_{\Delta_\ell})_{\ell\in N}$ is contained in the closed ball with radius $\sqrt{\psi(T)-\psi(0)}$ of the Hilbert space $\Ltwo(d\psi,[0,T],\R^n)$. As such, there exists a subsequence $N'$ of $N$ such that $(y_{\Delta_\ell})_{\ell\in N'}$ converges to $y$ weakly in $\Ltwo(d\psi,[0,T],\R^n)$. It then follows from Lemma~\ref{lem:xderpsiy} that
\beq\label{e: derv of x}
\dot{x}(t)=\dot{\psi}(t)y(t)
\eeq
for almost all $t\in\Gamma$. 

Now, let $t^*\in\Gamma$. Then, for every $\ell\in N$, there must exist $k_\ell\in\nset{K(\Delta_\ell)}$ with the property that $t_{k_{\ell}}<t^*< t_{k_{\ell}+1}$. Note that $\lim_{\ell\uparrow\infty} t_{k_{\ell}}=\lim_{\ell\uparrow\infty}t_{k_{\ell}+1}=t^*$ since $\abs{\Delta_\ell}$ converges to zero as $\ell$ tends to infinity. By construction, we have
$$
\big (x_{t_{k_{\ell}+1}}, -\frac{x_{t_{k_{\ell}+1}}-x_{t_{k_{\ell}}}}{t_{k_{\ell}+1}-t_{k_{\ell}}} \big )\in\graph F(t_{k_{\ell}+1},\cdot).
$$
Equivalently, we have
\begin{equation}\label{eq:test}
\big (x_{\Delta_\ell}(t_{k_{\ell}+1}), -\frac{\psi(t_{k_{\ell}+1})-\psi(t_{k_{\ell}})}{t_{k_{\ell}+1}-t_{k_{\ell}}}y_{\Delta_\ell}(t) \big )\in\graph F(t_{k_{\ell}+1},\cdot). 
\end{equation}

Let $S_\ell(t^*) := - \frac{t_{k_{\ell}+1}-t_{k_{\ell}}}{\psi(t_{k_{\ell}+1})-\psi(t_{k_{\ell}})} F(t_{k_{\ell}+1},x_{\Delta_\ell}(t_{k_{\ell}+1}))$. From \eqref{eq:test}, we have that $y_\ell(t^*) \in S_\ell(t^*)$. We now invoke  Lemma~9 and observe that $y(t^*) \in \cl \left( \conv \left(\limsup_{\ell \to \infty} S_\ell (t^*) \right)\right)$. Due to the outer-semicontinuity assumption, we have $\limsup_{\ell \to \infty} F(t_{k_{\ell}+1},x_{\Delta_\ell}(t_{k_{\ell}+1})) \subseteq F(t^*,x(t^*))$. The set $F(t^*,x(t^*))$ is closed and convex because of the maximal monotonicity property, and hence
\[
y(t^*) \in \frac{-1}{\dot \psi(t^*)} F(t^*,x(t^*)).
\]
Since $\dot \psi(t^*) \geq 1$, we get
\[
\dot{x}(t^*)\overset{\eqref{e: derv of x}}{=}\dot{\psi}(t^*)y(t^*) \in -F(t^*, x(t^*))
\]
for each $t^* \in \Gamma$.
\EP

\section{Extensions}\label{s:ext}
In this section, we extend the results of Theorem~\ref{thm:main}. First, we consider non-autonomous differential inclusions of the form
\begin{equation}\label{eq:incGen-nonautonomous}
\dot{x}(t)\in -F\big(t,x(t)\big)+u(t),\quad x(0)=x_0 
\end{equation}
where $F(t,\cdot):\R^n \rightrightarrows \R^n$ is maximal monotone for all $t\geq 0$ and $u \in \Lone([0,T],\R^n)$. We begin with the following observation.
\begin{lem}\label{l:extinput}
Consider the system \eqref{eq:incGen-nonautonomous}. Let $u \in \Lone([0,T],\R^n)$ and $G$ be the set-valued mapping defined by $G(t,\xi):=F\big(t,\xi+\int_0^t{u(\tau)\,d \tau}\big)$. Then, the differential inclusion \eqref{eq:incGen-nonautonomous}
admits a solution $x$ if and only if the differential inclusion
\begin{align}
	\dot \xi(t) \in - G\big(t,\xi(t)\big),\quad \xi(0)=x_0
\end{align}
admits a solution $\xi$.
\end{lem}
\BP For the `only if' part, suppose that $x$ is a solution of \eqref{eq:incGen-nonautonomous}. Define
$$
  \xi(t) = x(t) - \int_0^t{u(\tau)\,d \tau}
$$
for all $t\geq 0$. Note that $\xi(0)=x_0$ and
$$
  \dot \xi(t) = \dot x(t)-u(t) \in -F\big(t, x(t)\big)=-F\big(t,\xi(t)+\int_0^t{u(\tau)\,d\tau}\big) = - G(t,\xi(t)).
$$
The `if' part follows reversing the arguments.
\EP

\begin{thm}\label{t:non-auto}
Consider system~\eqref{eq:incGen-nonautonomous} and suppose that $u \in \Lone([0,T],\R^n)$ and $F(t,\cdot)$ satisfies assumptions \ref{ass:basicF}, \ref{ass:new}, and \ref{ass:dist}. If \ref{ass:least} holds, then for every $x_0\in\cl\big(\dom F(0,\cdot)\big)$, 
there exists a unique solution $x\in AC([0,T],\R^n)$.
If \ref{ass:static} holds, then for every $x_0\in \dom F(0,\cdot)$, 
there exists a unique solution $x\in AC([0,T],\R^n)$.
\end{thm}
\BP
In view of Lemma~\ref{l:extinput} and Theorem~\ref{thm:main}, it is enough to show that the time-dependent set-valued map $G$ defined by $G(t,x)=F\big(t,x+\Phi(t)\big)$ with $\Phi(t)=\int_0^t{u(\tau)\,d \tau}$ satisfies assumptions \ref{ass:basicF}, \ref{ass:new}, \ref{ass:dist}, and in addition, $G$ satisfies \ref{ass:least} (resp. \ref{ass:static}) if $F$ satisfies \ref{ass:least} (resp. \ref{ass:static}).\\

\ref{ass:basicF}: Since maximal monotonicity is preserved under translations (see e.g. \cite[Thm. 12.43]{Rockafellar98}), $G$ satisfies assumption \ref{ass:basicF} whenever $F$ satisfies it.\\ 

\ref{ass:new}: Note that $\graph G(t,\cdot)=\graph F(t,\cdot)-\pset{\Phi(t)}\times \pset{0}$. Since the set-valued mapping $t \mapsto \graph F(t,\cdot)$ is outer semicontinuous on $[0,T]$ by assumption and $\Phi$ is absolutely continuous, $t \mapsto \graph G(t,\cdot)$ is outer semicontinuous  on $[0,T]$.

\ref{ass:dist}: Note that $\dom G(t,\cdot)=\dom F(t,\cdot)-\Phi(t)$ for all $t\in[0,T]$. Therefore, we have
\begin{align*}
\sup_{z \in \dom G(s,\cdot)}\dist{\big(z, \dom G(t,\cdot)\big)}&=\sup_{z+\Phi(s) \in \dom F(s,\cdot)}\dist{\big(z+\Phi(t), \dom F(t,\cdot)\big)}\\
&\leq \varphi(t)-\varphi(s)+\abs{\Phi(t)-\Phi(s)}
\end{align*}
for all $s,\,t$ with $0\leq s\leq t \leq T$ since $F$ satisfies assumption \ref{ass:dist}. Note that
$$
\varphi(t)-\varphi(s)+\abs{\Phi(t)-\Phi(s)}\leq \bar \varphi(t)-\bar \varphi(s)
$$	
for all $s,\,t$ with $0\leq s\leq t \leq T$ where 
\beq\label{e:varphi from F to G}
\bar \varphi(t):=\varphi(t)+\int_0^t{\abs{u(\tau)}\,d \tau}.
\eeq
Therefore, $G$ satisfies assumption \ref{ass:dist}.\\

\ref{ass:least}: { Since $F$ satisfies assumption \ref{ass:least}, we observe that
\beq\label{e:assumption least}
\abs{G^0(t,x)}=\abs{F^0\big(t,x+\Phi(t)\big)}
\leq \sigma(t)\big(1+\abs{x+\Phi(t)}\big) 
\eeq
for all $x \in \dom{G(t,\cdot)}$ with $t\in[0,T]$.
Note that
$$
\sigma(t)\big(1+\abs{x+\Phi(t)}\big)\leq \sigma(t)\big(1+\abs{\Phi(t)}\big)(1+\abs{x}).
$$
Together with \eqref{e:assumption least}, this results in
\beq\label{e:least norm from F to G}
\abs{G^0(t,x)}\leq \sigma(t)\big(1+\abs{\Phi(t)}\big)(1+\abs{x})
\eeq
for all $t\in[0,T]$ and for all $x\in\dom G(t,\cdot)$. Since the function $t\mapsto \sigma(t)$ is integrable and $t\mapsto 1+\abs{\Phi(t)}$ is continuous, their product is integrable. Consequently, $G$ satisfies assumption \ref{ass:least}.

\ref{ass:static}: Let $Y_\lambda^F(t,\cdot)$ and $Y_\lambda^G(t,\cdot)$ denote the Yosida approximation of $F(t,\cdot)$ and $G(t,\cdot)$ respectively. With $F$ satisfying \ref{ass:static}, and observing that $x \in \dom G(s,\cdot)$ if and only if $x + \Phi(s) \in \dom F(s,\cdot)$, we get, for all $0 \le s < T$, $\lambda < \Lambda$,
\begin{align}
\abs{Y_\lambda^G(s+\lambda,x)} & \le \abs{Y_\lambda^G(s+\lambda,x+\Phi(s)-\Phi(s+\lambda))} + \abs{Y_\lambda^G(s+\lambda,x) - Y_\lambda^G(s+\lambda,x+\Phi(s)-\Phi(s+\lambda))} \notag \\
& \le \abs{Y_\lambda^F(s+\lambda,x+\Phi(s))} + M_u \notag \\
 & \le M_u + \theta (s+\lambda) - \theta(s) + \abs{F^0(s,x+\Phi(s))} \notag \\
& = \overline \theta (s+\lambda) - \overline \theta(s) + \abs{G^0(s,x)} \label{eq:defThetabarYB}
\end{align}
where $\overline \theta(s) = M_u + \theta(s)$ for some $M_u > 0$. In writing down the second inequality, we used the fact that $Y_\lambda^G(t,x) = Y_\lambda^F(t,x+\Phi(t))$, and moreover, $\lambda^{-1}$-Lipschitzian property of $Y_\lambda^G(t,\cdot)$ along with integrability of $u$ yields
\[
\abs{Y_\lambda^G(s+\lambda,x) - Y_\lambda^G(s+\lambda,x+\Phi(s)-\Phi(s+\lambda))} \le M_u
\]
for some $M_u > 0$.
\EP
}

Now, we turn our attention to differential inclusions of the form
\begin{equation}\label{eq:incGen-Lip}
\dot{x}(t)\in -F\big(t,x(t)\big)+f\big(x(t)\big)+u(t),\quad x(0)=x_0 
\end{equation}
where $F(t,\cdot):\R^n \rightrightarrows \R^n$ is maximal monotone for all $t\geq 0$, $f:\R^n\rightarrow\R^n$ is a function and $u \in \Lone([0,T],\R^n)$. Based on Theorem~\ref{t:non-auto}, we present the following existence and uniqueness result.


\begin{thm}
Consider system~\eqref{eq:incGen-Lip} and suppose that $u \in \Lone([0,T],\R^n)$, $f:\R^n\rightarrow\R^n$ is a Lipschitz continuous function, and $F(t,\cdot)$ satisfies assumptions \ref{ass:basicF}, \ref{ass:new}, and \ref{ass:dist}. If \ref{ass:least} holds, then for every $x_0\in\cl\big(\dom F(0,\cdot)\big)$, 
there exists a unique solution $x\in AC([0,T],\R^n)$.
If \ref{ass:static} holds, then for every $x_0\in \dom F(0,\cdot)$, 
there exists a unique solution $x\in AC([0,T],\R^n)$.
\end{thm}

\BP
Let $x_0\in\cl\big(\dom F(0,\cdot)\big)$ and let $x_0(t)=x_0$ for all $t\in[0,T]$. It follows from Theorem~\ref{t:non-auto} that for each integer $\ell$ with $\ell\geq 1$ there exists a unique absolutely continuous function $x_{\ell+1}:[0,T]\rightarrow\R^n$ such that $x_{\ell+1}(0)=x_0$, $x_{\ell+1}(t) \in\dom  F (t,\cdot)$ and the differential inclusion 
$$
\dot{x}_{\ell+1}(t)\in -F\big(t,x_{\ell+1}(t)\big)+f\big(x_{\ell}(t)\big)+u(t)
$$
holds for almost all $t \in[0,T]$. In the rest of the proof, we will construct a solution of \eqref{eq:incGen-Lip} by showing that the sequence $\{x_\ell(\tau)\}_{\ell \in \N}$ is a Cauchy sequence that converges to an absolutely continuous function which satisfies \eqref{eq:incGen-Lip}. 

{\em Step~1: The sequence $\{x_\ell(\tau)\}_{\ell \in \N}$ is Cauchy.\/}
By using \ref{ass:basicF} and Lipschitzness of $f$, we see that
\begin{align*}
\half\frac{d}{dt}\big(\abs{x_{\ell+1}(t)-x_{\ell}(t)}^2\big)&=\inn{\dot{x}_{\ell+1}(t)-\dot{x}_{\ell}(t)}{x_{\ell+1}(t)-x_{\ell}(t)}\\
&\leq L_f \abs{x_{\ell}(t)-x_{\ell-1}(t)}\,\abs{x_{\ell+1}(t)-x_{\ell}(t)}
\end{align*}
for almost all $t \in[0,T]$ where $L_f$ is the Lipschitz constant of $f$. By integrating both sides from $0$ to $\tau\in[0,T]$, we obtain
$$
\half \abs{x_{\ell+1}(\tau)-x_{\ell}(\tau)}^2\leq \int_0^\tau L_f \abs{x_{\ell}(s)-x_{\ell-1}(s)}\,\abs{x_{\ell+1}(s)-x_{\ell}(s)}\,ds.
$$
Application of \cite[Lemma A.5, p. 157]{Brez73} results in
$$
\abs{x_{\ell+1}(\tau)-x_{\ell}(\tau)}\leq \int_0^\tau L_f \abs{x_{\ell}(s)-x_{\ell-1}(s)}\,ds
$$
for all $\tau\in[0,T]$. Hence, we get
\beq\label{e: bound x l plus one to x l}
\abs{x_{\ell+1}(\tau)-x_{\ell}(\tau)}\leq \frac{(L_f \tau)^\ell}{\ell!}\abs{x_1-x_0}_{\Linf}
\eeq
for all $\tau\in[0,T]$ where $\abs{\cdot}_{\Linf}$ denotes the sup norm. Consequently, $(x_{\ell})_{\ell\in\N}$ converges uniformly on $[0,T]$ to a function $x$.

{\em Step~2: The function $x$ belongs to $AC([0,T],\R^n)$.\/} 
For the moment, suppose that there exist an integer $L$ and a nondecreasing function $\widehat\psi:AC([0,T],\R)$ such that
\beq\label{e:uni bound on x ell}
\abs{x_\ell(\ol \tau)-x_\ell(\ul\tau)}\leq \widehat\psi(\ol\tau) - \widehat\psi(\ul\tau)
\eeq
for all $\ell\geq L$ and for all $\ul\tau,\ol \tau$ with $0\leq \ul\tau < \ol \tau\leq T$.
This would mean that we have
\begin{align*}
\abs{x(\ol\tau)-x(\ul\tau)}&\leq\abs{x(\ol\tau)-x_{\ell}(\ol\tau)}+\abs{x_{\ell}(\ol\tau)-x_{\ell}(\ul\tau)}+\abs{x_{\ell}(\ul\tau)-x(\ul\tau)}\notag\\
&\leq\abs{x(\ol\tau)-x_{\ell}(\ol\tau)}+\widehat\psi(\ol\tau)-\widehat\psi(\ul\tau)+\abs{x_{\ell}(\ul\tau)-x(\ul\tau)}\notag\\
&\leq \widehat\psi(\ol\tau)-\widehat\psi(\ul\tau)
\end{align*} 
for all $\ul\tau,\ol \tau$ with $0\leq \ul\tau < \ol \tau\leq T$ where the first inequality follows from the triangle inequality, the second from \eqref{e:uni bound on x ell} for all $\ell\geq L$, and the third by taking the limit as $\ell$ tends to infinity. Thus, absolute continuity of $x$ follows from absolute continuity of the function $\widehat\psi$.

To prove \eqref{e:uni bound on x ell}, we first observe that the triangle inequality and \eqref{e: bound x l plus one to x l} result in 
\begin{align*}
\abs{x_\ell (\tau)} &\le \abs {x_\ell(\tau) - x_{\ell-1}(\tau)} + \abs{x_{\ell-1}(\tau)} \\
& \le \sum_{i=1}^{l-1} \abs {x_{i+1}(\tau) - x_{i}(\tau)} + \abs{x_1(\tau)} \\
& \le (e^{L_f \tau} -1) + \abs{x_1(\tau)} \le C_1
\end{align*}
for some $C_1 > 0$, and each $\ell \geq 1$, $\tau \in [0,T]$.
Hence, in particular, $\abs{x_\ell (\cdot)}$ is uniformly bounded for each $\ell \geq 1$.
The bound in~\eqref{e:uni bound on x ell} is obtained differently in case if $F(\cdot,\cdot)$ satisfies \ref{ass:least} or \ref{ass:static}.

{\em If $F$ satisfies \ref{ass:least}:} Similar to function $\ol \varphi$ in \eqref{e:varphi from F to G}, let us introduce the function $\varphi_\ell$ as,
\[
\varphi_\ell(t) = \varphi(t) + \int_0^t \abs{f(x_{\ell-1}(\tau))} + \abs{u(\tau)} \, d\tau 
\]
so that, the Lipschitz continuity of $f$, $\abs{f(z)} \le \abs{f(0)} + L_f \abs{z}$, yields
\begin{align*}
\varphi_\ell(\ol \tau) - \varphi_\ell(\ul \tau) & \le L_f \int_{\ul \tau}^{\ol \tau} \abs{x_{\ell-1}(\tau)}\,d\tau + \abs{f(0)} (\ol \tau - \ul \tau) + \int_{\ul \tau}^{\ol \tau} \abs{u(\tau)}\,d\tau \\
& \le (L_f C_1+ \abs{f(0)})(\ol \tau - \ul \tau) + \int_{\ul \tau}^{\ol \tau} \abs{u(\tau)}\,d\tau.
\end{align*}
By introducing the function $\psi_\ell$, similar to \eqref{eq:defPsi}, as
\[
\psi_\ell(t) = t+2\varphi_\ell(t)+ (1+\gamma)\int_0^t \sigma(s)\,ds\quad\forall\,t\in[0,T],
\]
and by letting $g(s):= s + (1+\gamma)\int_0^s \sigma(\tau)\,d\tau + 2 \int_0^s \abs{u(\tau)}\,d\tau$, we get 
\begin{align*}
\psi_\ell( \ol \tau) - \psi_\ell (\ul \tau) &= \ol \tau - \ul \tau + 2 \varphi_\ell(\ol \tau) - 2\varphi_\ell(\ul \tau) + (1+\gamma)\int_{\ul \tau}^{\ol \tau} \sigma(\tau) d\tau \\
& \le g(\ol \tau) - g(\ul \tau) + 2 ( L_f C_1 + \abs{f(0)})(\ol \tau - \ul \tau) .
\end{align*}
It follows from \eqref{e:bound x with psi} that $\abs{x_\ell( \ol \tau) - x_\ell (\ul \tau)} \le \psi_\ell( \ol \tau) - \psi_\ell (\ul \tau)$. Thus we get
\[
\abs{x_\ell( \ol \tau) - x_\ell (\ul \tau)} \le g(\ol \tau) - g(\ul \tau) + 2 \, L_f C_1(\ol \tau - \ul \tau)
\]
for each $\ell \geq 1$, and hence \eqref{e:uni bound on x ell} follows with $\widehat \psi (s) := g(s) + 2\,L_fC_1 s$, which is clearly an absolutely continuous function.

{{\em If $F$ satisfies \ref{ass:static}:} In that case, the function $\psi_\ell$ is obtained from \eqref{eq:defPsiYB}, so that $\psi_{\ell}(s) := \overline\alpha \, s$. The constant $\overline \alpha:= 2 (\overline \theta(T) - \overline \theta(0))$, where $\overline\theta$ is defined as in \eqref{eq:defThetabarYB}, that is, $\overline \theta(s) := M_u + C_1 + \theta(s)$, for some $C_1, M_u \ge 0$, which is independent of $\ell$. Hence, we take $\widehat \psi(s) = \overline \alpha\, s$, which is clearly absolutely continuous.}

{\em Step~3: The function $x$ satisfies \eqref{eq:incGen-Lip}.}
In view of Lemma~\ref{l:extinput}, $x_\ell$ is a solution to the differential inclusion
$$
\dot{x}_{\ell}(t)\in -G_\ell\big(t,x_\ell(t)\big),\quad x_\ell(0)=x_0
$$
for all $\ell\geq 1$ where $G_\ell(t,\xi)=F(t,\xi+\Phi_\ell(t)\big)$ with $\Phi_\ell(t)=\int_0^t f\big(x_{\ell-1}(\tau)\big)+u(\tau)\,d\tau$. Since $\graph G_\ell(t,\cdot)=\graph F(t,\cdot)-\pset{\Phi_\ell(t)}\times\pset{0}$.

Let us introduce the sequence $y_\ell := \frac{d x_\ell}{d\widehat \psi}$, so that
\[
\dot x_\ell (t) = \dot {\widehat \psi} (t) y_\ell (t) \in -F(t,x_\ell(t)) + f(x_\ell(t)) + u(t).
\]
Because of the bound \eqref{e:uni bound on x ell},  $\abs{y_\ell}_{L^\infty} \le 1$, and there exists a subsequence $N$ such that $(y_{\ell})_{\ell\in N}$ converges to $y$ weakly in $\Ltwo(d\widehat \psi,[0,T],\R^n)$.  It then follows from Lemma~\ref{lem:xderpsiy} that
\beq
\dot{x}(t)=\dot{\widehat \psi}(t)y(t)
\eeq
for almost all $t\in\Gamma :=\{t\in [0,T] \, \vert\, x_\ell, \ell \geq L, x, \text{ and } \psi \text{ are differentiable at } t\}$.

Let $t^* \in \Gamma$. By construction, we have
\[
\left( x_\ell(t^*), - \dot x_\ell (t^*) \right) \in \graph G_\ell \big(t^*, x_\ell(t^*)\big)
\]
or equivalently, 
\[
\left( x_\ell(t^*), - \dot {\widehat \psi}(t^*) y_\ell (t^*) \right) \in \graph G_\ell (t^*, x_\ell(t^*)).
\]
In other words, $y_\ell (t^*)$ belongs to the convex set $\frac{-1}{\dot {\widehat \psi}(t^*)}G_\ell (t^*, x_\ell(t^*))$. Using Lemma~\ref{l:values of weak limit-general} with $S_\ell(t) = \frac{-1}{\dot {\widehat \psi}(t)}G_\ell (t, x_\ell(t))$, and recalling that $\abs{\dot{\widehat \psi}(t)}\geq 1$ for each $t\in [0,T]$, we see that
\[
\dot x(t^*) = \dot {\widehat \psi}(t^*) y (t^*) \in \cl\left(\conv\left(\limsup_{\ell\rightarrow\infty}S_\ell(t^*)\right)\right) \subseteq - G(t^*, x(t^*))
\]
and the same holds for almost every $t^* \in [0,T]$.
\EP

\section{Linear systems and maximal monotone relations}\label{sec:compSys}

A particularly interesting class of time-dependent maximal monotone mappings arises from the interconnection of linear passive systems with maximal monotone relations. To formalize this class of systems, consider the linear system
\bse\label{e:lin+rel}
\begin{gather}
\dot{x}(t)=Ax(t)+Bz(t)+u(t)\label{e:lin+rel.1}\\
w(t)=Cx(t)+Dz(t)+v(t)\label{e:lin+rel.2}
\end{gather}
where $x\in\R^n$ is the state, $u\in\R^{n}$ and $v\in\R^m$ are external inputs, and $(z,w)\in\R^{m+ m}$ are the external variables that satisfy
\beq\label{e:lin+rel-td.3}
\big(-z(t),w(t)\big)\in \graph(\cM)
\eeq
\ese
for some set-valued map $\cM:\R^m\rightrightarrows\R^m$.

By solving $z$ from the relations \eqref{e:lin+rel.2}, \eqref{e:lin+rel-td.3}, and substituting in \eqref{e:lin+rel.1}, we obtain the differential inclusion
\beq\label{e:dispec}
\dot{x}(t)\in -H\big(t,x(t)\big)+u(t)
\eeq
where
\beq\label{e:defht}
H(t,x)=-Ax+B(\cM+D)\inv\big(Cx+v(t)\big)
\eeq
and
\beq\label{e:dom-ht}
\dom H(t,\cdot)=C\inv\big(\im(\cM+D)-v(t)\big).
\eeq

The rest of this section is devoted to developing conditions under which the time-dependent set-valued mapping $H(t,\cdot)$ satisfies the hypotheses of Theorem~\ref{thm:main}. To establish such conditions, we first introduce the notion of {\em passivity} of a linear system.

A linear system $\Sigma(A,B,C,D)$
\bse\label{e:lin}
\begin{gather}
\dot{x}(t)=Ax(t)+Bz(t)\label{e:lin1}\\
w(t)=Cx(t)+Dz(t)\label{e:lin2}
\end{gather}
\ese
is said to be {\em passive}, if there exists a nonnegative-valued {\em storage function\/} $V:\R^n\rightarrow\R_+$ such that the so-called {\em dissipation inequality\/}
\beq
V(x(t_1))+\int_{t_1}^{t_2}z^\top (\tau)w(\tau)\,d\tau\leq V(x(t_2))
\eeq
holds for all $t_1,t_2$ with $t_1< t_2$ and for all trajectories $(z,x,w)\in\Ltwo([t_1,t_2],\R^m)\times AC([t_1,t_2],\R^n)\times\Ltwo([t_1,t_2],\R^m)$ of the system \eqref{e:lin}.

The classical Kalman-Yakubovich-Popov lemma states that the system \eqref{e:lin} is passive if, and only if, the linear matrix inequalities
\beq\label{e:lmis}
K=K^\top \geq 0\qquad\bbm A^\top K+KA & KB-C^\top \\B^\top K-C& -(D^\top +D)\ebm\leq 0
\eeq
admits a solution $K$. Moreover, $V(x)=\half x^\top Kx$ defines a storage function in case $K$ is a solution of the linear matrix inequalities \eqref{e:lmis}.

In the following proposition, we summarize some of the consequences of passivity that will be used later.

\begin{proposition}[{\cite[Lem. 1]{camlibel:14}}]\label{p:pass}
If $\Sigma(A,B,C,D)$ is passive with the storage function $x\mapsto\half x^\top Kx$ then the following statements hold:
\begin{enumerate}
\renewcommand{\theenumi}{\roman{enumi}.}
\renewcommand{\labelenumi}{\theenumi}
\item $D$ is positive semi-definite.
\item $(KB-C^\top )\ker(D+D^\top )=\{0\}$.
\end{enumerate}
\end{proposition}

The following theorem states conditions that guarantee the hypotheses of Theorem~\ref{t:non-auto} for the time-dependent set-valued mapping $H$ as defined in \eqref{e:defht}.

\begin{thm}\label{thm: linear case}
Let $T>0$. Suppose that
\begin{enumerate}
\renewcommand{\theenumi}{\roman{enumi}.}
\renewcommand{\labelenumi}{\theenumi}
\item\label{i:passivity} $\Sigma\abcd$ is passive with the storage function $x\mapsto \half x^\top x$,
\item\label{i: max mon} $\cM$ is maximal monotone,
\item\label{i: rel int} for all $t\in[0,T]$, we have\footnote{Here, $\rint(S)$ denotes the relative interior of $S$.} $\im C\cap\rint\big(\im(\cM+D)-v(t)\big)\neq \emptyset$, 
\item\label{i: v(t)} $v$ is bounded on $[0,T]$,
\item\label{i: dom abs cont} there exists an absolutely continuous nondecreasing function $\theta : [0,T] \rightarrow \R$ such that
$$
\sup_{w\in \im C\cap \left(\im(\cM+D)-v(s)\right)}\dist\Big(w,\im C\cap \big(\im(\cM+D)-v(t)\big)\Big)\leq \theta(t)- \theta(s)
$$
for all $s,t$ with $0\leq s \leq t\leq T$.
\item\label{i: H least} { There exists a positive number $\alpha$ such that 
	\begin{equation}\label{eq:least-spec}
	\abs{B\big((\calM+D)\inv\big)^0(\eta)}\leq \alpha (1+\abs{\eta})
  \end{equation}
for all $\eta \in \dom{(\calM+D)\inv}$.
}
\end{enumerate}
Then, $H$ satisfies assumptions \ref{ass:basicF}, \ref{ass:new}, \ref{ass:dist}, and \ref{ass:least}.
\end{thm}

\BP For brevity, we define $W(t):=\im(\cM+D)-v(t)$ for all $t\in[0,T]$. Note that $\dom H(t,\cdot)=C\inv W(t)$ for all $t\in[0,T]$.\\

\ref{ass:basicF}: It follows from \cite[Thm. 2]{CamlSchu16} that the conditions \eqref{i:passivity}-\eqref{i: rel int} imply that $H(t,\cdot)$ is a maximal monotone mapping for all $t\in[0,T]$. As such, $H$ satisfies assumption \ref{ass:basicF}.\\

\ref{ass:dist}: Let $t$ and $s$ be such that $0\leq s\leq t\leq T$. Also, let $x\in C\inv W(s)$ and let $y=\proj\big(x,C\inv W(t)\big)$. Further, let $\zeta=\proj\big(Cx,\im C\cap W(t)\big)$. Therefore, there exists $\xi$ such that $\zeta=C\xi$. Without loss of generality, we can assume that $x-\xi\in\im C^\top $ since $\R^n=\im C^\top \oplus \ker C$. Now, we see that $Cx\in\im C\cap W(s)$ and $\zeta=C\xi\in\cl\big(\im C\cap W(t)\big)$. From \eqref{i: dom abs cont}, we get
\beq\label{e: bound on Cx and Cxi}
\abs{Cx-C\xi}\leq \theta(t)- \theta(s).
\eeq
Since $x-\xi\in\im C^\top $, there exists a positive number $\alpha$ such that
\beq
\abs{x-\xi}\leq \alpha \big( \theta(t)- \theta(s)\big).
\eeq
Since $\xi\in C\inv W(t)$, we obtain
$$
\abs{x-y}\leq\abs{x-\xi}\leq \alpha \big( \theta(t)- \theta(s)\big).
$$
Therefore, we see that
$$
\dist\big(x,C\inv W(t)\big)\leq \alpha \big( \theta(t)- \theta(s)\big).
$$
This implies that
$$
\sup_{x\in C\inv W(s)}\dist\big(x,C\inv W(t)\big)\leq \alpha \big( \theta(t)- \theta(s)\big).
$$
Since $\dom H(t,\cdot)=C\inv W(t)$, we can conclude that $H$ satisfies assumption \ref{ass:dist}.\\

\ref{ass:least}: {  From \eqref{i: H least}, we know that there exists a positive number $\alpha$ such that 
	\begin{equation}\label{eq:least-spec-again}
	\abs{B\big((\calM+D)\inv\big)^0(\eta)}\leq \alpha(1+\abs{\eta})
  \end{equation}
for all $\eta \in \dom{(\calM+D)\inv}$. Let $x \in \dom{H(t,\cdot)}$ for some $t\in[0,T]$. Since $Ax-Bz^0\in H(t,x)$ where $z^0\in\big((\calM+D)\inv)^0\big(Cx+v(t)\big)$, we have
\beq\label{e:hopefully last label}
\abs{H^0(t,x)}\leq \abs{Ax-Bz^0}\leq \abs{Ax}+\abs{Bz^0}.
\eeq
Since $x \in \dom{H(t,\cdot)}$, it follows that $Cx+v(t)\in\mathbb{B}^m(\rho)\cap \dom(\calM+D)\inv$,  and using \eqref{eq:least-spec-again}, \eqref{e:hopefully last label}, and boundedness of $v$, we get
$$
\abs{H^0(t,x)}\leq \abs{Ax}+\alpha_\rho(1+\abs{Cx+v(t)})\leq \beta (1+\abs{x})
$$
for some positive number $\beta$ that does not depend on $t$. Therefore, we have
$$
\abs{H^0(t,x)}\leq \beta (1+\abs{x})
$$
for all $x \in \dom{H(t,\cdot)}$ with $t\in[0,T]$. In other words, $H$ satisfies assumption  \ref{ass:least}.\\
}

\ref{ass:new}: Let $(t_\ell,x_\ell,y_\ell)_{\ell\in \N}\subseteq[0,T]\times\R^n\times\R^n$ be a sequence such that $y_\ell\in H(t_\ell,x_\ell)$ and $\lim_{\ell\uparrow\infty}(t_\ell,x_\ell,y_\ell)=(t,x,y)$ for some $t\in[0,T]$, $x,y\in\R^n$. What needs to be proven is that $y\in H(t,x)$.

Note that for each $\ell$ there exists $z_{\ell}\in(\cM+D)\inv\big(Cx_{\ell}+v(t_\ell)\big)$ such that $y_{\ell}=-Ax_{\ell}+Bz_{\ell}$. Then, $(Bz_{\ell})_{\ell\in\N}$ converges. Let $\calW$ be the subspace parallel to the affine hull of $\im(\cM+D)=\dom(\cM+D)\inv$. It follows from maximal monotonicity of $(\cM+D)\inv$ that for each $\ell$
\beq\label{e:dommd}
\zeta+z_{\ell}\in(\cM+D)\inv\big(Cx_{\ell}+v(t_\ell)\big)
\eeq
holds for any $\zeta\in\calW^\perp$. Now, let $z_{\ell}=z_{\ell}^1+z_{\ell}^2$ where $z_{\ell}^1\in\ker B\cap\calW^\perp$ and
\beq\label{e:z2 space}
z_{\ell}^2\in(\ker B\cap\calW^\perp)^\perp=\im B^\top +\calW.
\eeq
Note that
\beq\label{e:bznu=bz2nu}
Bz_{\ell}=Bz_{\ell}^2.
\eeq
From \eqref{e:dommd} by taking $\zeta=-z_{\ell}^1$, we have
\beq\label{e: z 2 ell}
z_{\ell}^2\in(\cM+D)\inv(Cx_{\ell}+v(t_\ell)).
\eeq
Suppose that $(z^2_{\ell})_{\ell\in\N}$ is bounded. Hence, $(z^2_{\ell})_{\ell\in\N}$ converges on a subsequence $N$, say to $z$. Then, we see from $(x_\ell,y_\ell)=(x_\ell,-Ax_\ell+Bz_\ell)$ that $(x_\ell,y_\ell)_{\ell\in N}$ converges to $(x,y)=(x,-Ax+Bz)$. Since $H(t,\cdot)$ is maximal monotone and that implies the closedness of $\graph\big(H(t,\cdot)\big)$, we then can conclude that $(x,y)\in\graph\big(H(t,\cdot)\big)$, or equivalently $y\in H(t,x)$. 

Therefore, it is enough to show that $(z^2_{\ell})_{\ell\in\N}$ is bounded. Suppose, on the contrary, that $z^2_{\ell}$ is unbounded. Without loss of generality, we can assume that the sequence $\frac{z^2_{\ell}}{\abs{z^2_{\ell}}}$ converges. Define
\beq\label{e:zinfnorm1}
\zeta_\infty=\lim_{\ell\rightarrow\infty} \frac{z^2_{\ell}}{\abs{z^2_{\ell}}}.
\eeq
From \eqref{e:bznu=bz2nu} and the fact that $(Bz_{\ell})_{\ell\in\N}$ converges, we have
\beq\label{e:closedness-1}
\lim_{\ell\rightarrow\infty}Bz^2_{\ell}=B\zeta.
\eeq
Thus, we get
\beq\label{e:55}
\zeta_\infty\in\ker B.
\eeq
Let $(\barx,\bary)\in\graph H(t,\cdot)$. Then, $\bary=-A\barx+B\barz$ where
\beq\label{e: z of y}
\barz\in (\cM+D)\inv \big(C\barx+v(t)\big).
\eeq

Due to passivity with $K=I$ and monotonicity of $(\cM+D)\inv$, it follows from \eqref{e: z 2 ell} and \eqref{e: z of y} that
\begin{align*}
\inn{x_{\ell}-\barx}{-A(x_{\ell}-\barx)+B(z^2_{\ell}-\barz)}&\geq\inn{z^2_{\ell}-\barz}{C(x_{\ell}-\barx)-D(z^2_{\ell}-\barz)}\\
&\geq -\inn{z^2_{\ell}-\barz}{v(t_\ell)-v(t)}.
\end{align*}
By dividing by $\abs{z^2_{\ell}}^2$, taking the limit as $\ell$ tends to infinity and using boundedness of $v$, we obtain
\beq
\inn{\zeta_\infty}{D\zeta_\infty}\leq 0.
\eeq
Since $D$ is positive semi-definite due to the first statement of Proposition~\ref{p:pass}, this results in
\beq
\zeta_\infty\in\ker(D+D^\top ).
\eeq
Then, it follows from \eqref{e:55}, $K=I$, and the second statement of Proposition~\ref{p:pass} that
\beq\label{e:zinfkerct}
\zeta_\infty\in\ker C^\top .
\eeq
Let $\eta\in\im(\cM+D)-v(t)$ and $\zeta\in(\cM+D)\inv\big(\eta+v(t)\big)$. In view of monotonicity of $(\cM+D)\inv$, we get
\beq
\inn{\frac{z^2_{\ell}-\zeta}{\abs{z^2_{\ell}}}}{Cx_{\ell}+v(t_\ell)-\eta-v(t)}\geq 0,
\eeq
from \eqref{e: z 2 ell}. Taking the limit as $\ell$ tends to infinity, employing boundedness of $v$, and using \eqref{e:zinfkerct}, we obtain
\beq\label{e:zeta inf separates}
\inn{\zeta_\infty}{Cx-\eta}=\inn{\zeta_\infty}{-\eta}\geq 0.
\eeq
From \eqref{e:zinfkerct} and \eqref{e:zeta inf separates}, we see that the hyperplane $\spn(\pset{\zeta_\infty})^\perp$ separates the sets $\im C$ and $\im(\cM+D)-v(t)$. In view of $\im C=\rint(\im C)$ and \eqref{i: rel int}, it follows from \cite[Thm. 11.3]{rockafellar:70} that $\im C$ and $\im(\cM+D)-v(t)$ cannot be properly separated. Therefore, both $\im C$ and $\im(\cM+D)-v(t)$ must be contained in the hyperplane $\spn(\pset{\zeta_\infty})^\perp$. Thus, we see that $\im(\cM+D)$ is contained in $v(t)+\spn(\pset{\zeta_\infty})^\perp$. Since $\calW$ is the subspace parallel to the affine hull of $\im(\cM+D)$, we get $\calW\subseteq\spn(\pset{\zeta_\infty})^\perp$ which implies $\zeta_\infty\in\calW^\perp$. Together with \eqref{e:55}, we get
$$
\zeta_\infty\in\ker B\cap\calW^\perp.
$$
In view of \eqref{e:z2 space} and \eqref{e:zinfnorm1}, this yields $\zeta_\infty=0$. This, however, clearly contradicts with \eqref{e:zinfnorm1}. Therefore, $\abs{z^2_{\ell}}$ must be bounded.\EP

Next, we specialize the results of Theorem~\ref{thm: linear case} to linear complementarity systems.
\subsection{Linear complementarity systems}
Linear complementarity systems are important instances of the differential inclusions of the form \eqref{e:dispec} with $\cM$ described by so-called complementarity relations. In this section, we aim at presenting tailor-made conditions for existence and uniqueness of solutions to linear complementarity systems. 

Consider a linear complementary system
\bse\label{e:lin+comp}
\begin{gather}
\dot{x}(t)=Ax(t)+Bz(t)+u(t)\label{e:lin+comp.1}\\
w(t)=Cx(t)+Dz(t)+v(t)\label{e:lin+comp.2}
\end{gather}
where $x\in\R^n$ is the state, $u\in\R^{n}$ and $v\in\R^m$ are external inputs, and $(z,w)\in\R^{m+ m}$ are the external variables that satisfy
\beq\label{e:lin+comp.3}
\big(-z(t),w(t)\big)\in \graph(\calP)
\eeq
\ese
where $\calP:\R^m\rightrightarrows\R^m$ is the maximal monotone set-valued mapping given by 
$$
\calP(\zeta)=\set{\eta}{\eta\geq 0,\,\zeta\leq 0,\text{ and }\inn{\eta}{\zeta}=0}.
$$
Next, we introduce the linear complementarity problem.

Given a vector $q\in\R^m$ and a matrix $M\in\R^{m\times m}$, the linear complementarity problem LCP($q,M$) is to find a vector $z\in\R^m$ such that
\begin{subequations}\label{e:lcp}
\begin{gather}
z\geq 0\label{e:cp1}\\
q+Mz\geq 0\label{e:cp2}\\
\inn{z}{q+Mz}=0.\label{e:cp3}
\end{gather}
\end{subequations}
We say that the LCP($q,M$) is {\em feasible\/} if there exists $z$ satisfying \eqref{e:cp1} and \eqref{e:cp2}. If a vector $z$ is feasible and satisfies \eqref{e:cp3} in addition, then we say that $z$ {\em solves\/} (is a {\em solution of\/}) LCP($q,M$). The set of all solutions of LCP($q,M$) will be denoted by SOL($q,M$).

A comprehensive study on LCPs can be found in \cite{cottle:92}. In the sequel, we will be interested in LCP($q,M$) where $M$ is a (not necessarily symmetric) positive semi-definite matrix.

Given a square matrix $M$, we define
\beq\label{e:qm-general} 
\calQ_M:=\mathrm{SOL}(0,M)=\set{z}{z\geq 0,\,Mz\geq 0,\text{ and }\inn{z}{Mz}=0} 
\eeq
and its dual cone
\beq\label{e:qm-dual} 
\calQ_M^+=\set{\zeta}{\inn{\zeta}{z}\geq 0\text{ for all }z\in\calQ_M}. 
\eeq
When $M$ is (not necessarily symmetric) a positive semi-definite matrix, the set $\calQ_M$ is a convex cone and can be given by $\calQ_M=\set{z}{z\geq 0,\,Mz\geq 0,\text{ and }(M+M^\top )z=0}$.

The following proposition characterizes the conditions under which an LCP with positive semi-definite $M$ matrix has solutions.

\begin{prop}[Cor. 3.8.10 of \cite{cottle:92} and Lem. 23 of \cite{camlibel:02a}]\label{p:lcp2} Let $M$ be a positive semi-definite matrix. Then, the following statements are equivalent:
\begin{enumerate}
\renewcommand{\theenumi}{\roman{enumi}.}
\renewcommand{\labelenumi}{\theenumi}
\item\label{lcp1} $q\in\calQ_M^+$.
\item\label{lcp2} $\mathrm{LCP}(q,M) \text{ is feasible}$.
\item\label{lcp3} $\mathrm{LCP}(q,M) \text{ is solvable}$.
\end{enumerate}
Moreover, the following statements hold:
\begin{enumerate}
\renewcommand{\theenumi}{\roman{enumi}.}
\renewcommand{\labelenumi}{\theenumi}
\setcounter{enumi}{3}
\item For each $q\in\calQ_M^+$, there exists a unique least-norm solution $z^*(q)\in\mathrm{SOL}(q,M)$ in the sense that $\abs{z^*(q)}\leq\abs{z}$ for all $z\in\mathrm{SOL}(q,M)$.
\item There exists a positive number $\alpha$ such that
$$
\abs{z^*(q)}\leq \alpha \abs{q}\qquad \forall\,q\in\calQ_M^+.
$$
\end{enumerate}
\end{prop}

Now, define 
\beq\label{e:def-comp}
H_\calP(t,x)=-Ax+B(\calP+D)\inv\big(Cx+v(t)\big).
\eeq
Note that $\dom H_\calP(t,\cdot)=C\inv\big(\im(\calP+D)-v(t)\big)$. Moreover, $q\in(\calP+D)(z)$ if and only if $-z\in\mathrm{SOL}(q,D)$. This means that $q\in(\calP+D)(z)$ if and only if $q\in\calQ_D^+$ in view of Proposition~\ref{p:lcp2}. In other words, $\dom H_\calP(t,\cdot)=C\inv\big(\calQ_D^+-v(t)\big)$.

The following theorem provides streamlined conditions that guarantee the hypotheses of Theorem~\ref{thm: linear case} for the time-dependent set-valued mapping $H_\calP$ as defined in \eqref{e:def-comp}.

\begin{thm}\label{e: linear case-comp}
Let $T>0$. Suppose that
\begin{enumerate}
\renewcommand{\theenumi}{\roman{enumi}.}
\renewcommand{\labelenumi}{\theenumi}
\item\label{i:passivity-2} $\Sigma\abcd$ is passive with the storage function $x\mapsto \half x^\top x$,
\item\label{i: rel int-2} $\im C\cap\rint\big(\im(\calP+D)-v(t)\big)\neq \emptyset$ for all $t\in[0,T]$,
\item\label{i: v(t)-2} $v\in AC([0,T],\R^m)$,
\end{enumerate}
Then, $H_\calP$ satisfies assumptions \ref{ass:basicF}, \ref{ass:new}, \ref{ass:dist}, and \ref{ass:least}.
\end{thm}

\BP It is enough to show that $H_\calP$ satisfies the hypotheses \eqref{i:passivity}-\eqref{i: H least} of Theorem~\ref{thm: linear case}. The first four hypotheses of Theorem~\ref{thm: linear case} are readily satisfied. Therefore, we need to show that the remaining two also hold.

For the hypothesis \eqref{i: dom abs cont} of Theorem~\ref{thm: linear case}, we need a streamlined version of Hoffman's bound on the polyhedral sets. To elaborate, let $\emptyset\neq\calR\subseteq\R^m$ be a polyhedral set given by $\calR=\set{\zeta}{R\zeta=0\text{ and }Q\zeta\leq q}$ where $R,Q$ are matrices and $q$ is vector with appropriate sizes. Hoffman's bound (see e.g. \cite[Lemma 3.2.3]{pang}) asserts that there exists a positive number $\alpha$ that depend on $\calR$ such that
\beq\label{e:hoffman}
\dist(x,\calR)\leq \alpha \big(\abs{Rx}+\abs{\max(0,Qx-q)}\big)
\eeq
for all $x\in\R^m$ where $\max$ denotes componentwise maximum. By definition $\calQ_D^+$ is a polyhedral cone. Therefore, we have $\calQ_D^+=\set{\eta\in\R^m}{Q\eta\leq 0}$ for some matrix $Q$. Let $E$ be a matrix such that $\im C=\ker E$. Then, we have
\begin{align}\label{e:poly domain}
\im C\cap \big(\im(\calP+D)-v(t)\big)&=\im C\cap \big(\calQ_D^+-v(t)\big)\\
&=\set{\zeta\in\R^m}{E\zeta=0\text{ and }Q\zeta\leq -Qv(t)}\notag
\end{align}
for all $t\in[0,T]$.

Let $s,t$ be such that $0\leq s\leq t\leq T$ and $w\in \im C\cap \big(\im(\calP+D)-v(s)\big)$. From \eqref{e:poly domain}, we see that
\beq\label{e:ineq for w}
Ew=0\text{ and }Qw\leq -Qv(s).
\eeq
Now, we have
\begin{align}
\dist\Big(w,\im C\cap \big(\im(\calP+D)-v(t)\big)\Big)&\overset{\eqref{e:hoffman}}{\leq}\alpha \Big(\abs{Ew}+\abs{\max\big(0,Qw+Qv(t)\big)}\Big)\nonumber\\
&\overset{\eqref{e:ineq for w}}{\leq}\alpha\abs{\max\big(0,-Qv(s)+Qv(t)\big)}\nonumber\\
&\leq\beta \abs{v(s)-v(t)}\label{e:almost there}
\end{align}
where $\beta$ is a positive number. Since $v$ is absolutely continuous, we have
$$
\abs{v(s)-v(t)}=\left|\int_s^t\dot{v}(\tau)\,d\tau\right|\leq \int_s^t\abs{\dot{v}(\tau)}\,d\tau=\int_0^t\abs{\dot{v}(\tau)}\,d\tau-\int_0^s\abs{\dot{v}(\tau)}\,d\tau.
$$
Then, \eqref{e:almost there} implies that
$$
\sup_{w\in \im C\cap \left(\im(\calP+D)-v(s)\right)}\dist\Big(w,\im C\cap \big(\im(\calP+D)-v(t)\big)\Big)\leq \theta(t)- \theta(s)
$$
for all $s,t$ with $0\leq s\leq t\leq T$ where $\theta(t)=\frac{1}{\beta}\int_0^t\abs{\dot{v}(\tau)}\,d\tau$ for all $t\in[0,T]$. Clearly, $\theta$ is nondecreasing and absolutely continuous.\\

For the hypothesis \eqref{i: H least} of Theorem~\ref{thm: linear case}, note that $\zeta\in(\calP+D)\inv(\eta)$ if and only if $-\zeta\in\mathrm{SOL}(\eta,D)$. Due to Proposition~\ref{p:pass}, $D$ is positive semi-definite. Then, it follows from Proposition~\ref{p:lcp2} that there exists a positive number $\alpha$ such that
$$
\big((\calP+D)\inv\big)^0(\eta)\leq \alpha \abs{\eta}
$$
for all $\eta\in\dom{(\calP+D)\inv}$. Therefore, $H_\calP$ satisfies the hypothesis \eqref{i: H least} of Theorem~\ref{thm: linear case}.
\EP

\section{Conclusions}\label{s:conc}

In this article, we have studied the existence of solutions to differential inclusions with time-dependent maximal monotone operators. With the help of an example, it is shown that our proposed conditions overcome the limitations of existing results. As a particular class of these inclusions, we consider differential equations coupled with time-dependent complementarity relations. For this system class, conditions for existence of solutions are derived explicitly in terms of system data.
To build on these results, the conditions for existence of solutions can be relaxed for differential inclusions where the maximal monotone operators have a particular structure, for example \cite{TanwBrog18}. 

Moving forward from the question of existence of solutions, it is also of interest to study the qualitative properties of the solutions of such systems, such as continuity with respect to initial data \cite{pang2009}. One can also investigate stability related problems for the generic class of dynamical systems, as has been done for some specific set-valued systems in \cite{TanwBrog14}. It also remains to be seen whether our proposed results provide any advantages in the study of optimal control problems, such as \cite{Briceno-Arias15}.



\bibliographystyle{plain}
\bibliography{monotone}

\begin{thebibliography}{10}

\bibitem{addi11}
K.~Addi, B.~Brogliato, and D.~Goeleven.
\newblock A qualitative mathematical analysis of a class of linear variational
  inequalities via semi-complementarity problems: {A}pplications in
  electronics.
\newblock {\em Mathematical Programming}, 126(1):31--67, 2011.

\bibitem{Adly18}
S.~Adly.
\newblock {\em A Variational Approach to Nonsmooth Dynamics: {A}pplications in
  Unilateral Mechanics and Electronics}.
\newblock SpringerBriefs in Mathematics. Springer International Publishing,
  2018.

\bibitem{AdlyHadd18}
S.~Adly, T.~Haddad, and B.-K. Le.
\newblock State-dependent implicit sweeping process in the framework of
  quasistatic evolution quasi-variational inequalities.
\newblock {\em Journal of Optimization Theory and Applications}, 2018.

\bibitem{AdlyHadd14}
S.~Adly, T.~Haddad, and L.~Thibault.
\newblock Convex sweeping process in the framework of measure differential
  inclusions and evolution variational inequalities.
\newblock {\em Mathematical Programming, Ser.~B}, 148(1):5--47, 2014.

\bibitem{Arseni-Benou99}
K.~Arseni-Benou, N.~Halidias, and N.S. Papageorgiou.
\newblock Nonconvex evolution inclusions generated by time-dependent
  subdifferential operators.
\newblock {\em Journal of Applied Mathematics and Stochastic Analysis},
  12(3):233--252, 1999.

\bibitem{attouch2019}
H.~Attouch and J.~Peypouquet.
\newblock Convergence of inertial dynamics and proximal algorithms governed by
  maximally monotone operators.
\newblock {\em Mathematical Programming, Ser.~B}, 174(1-2):391--432, 2019.

\bibitem{BausComb17}
H.H. Bauschke and P.L. Combettes.
\newblock {\em Convex Analysis and Monotone Operator Theory in {Hilbert}
  Spaces}.
\newblock Springer, NewYork, 2nd edition, 2017.

\bibitem{bertsekas}
D.P. Bertsekas.
\newblock {\em Convex Optimization Theory}.
\newblock Athena Scientific, 2009.

\bibitem{Brez73}
H.~Br\'ezis.
\newblock {\em Op\'erateurs Maximaux Monotones et Semi-Groupes des Contractions
  dans les Espaces de Hilbert}.
\newblock North-Holland, Mathematics Studies, 1973.

\bibitem{Briceno-Arias15}
L.M. Briceno-Arias, N.D. Hoang, and J.~Peypouquet.
\newblock Existence, stability and optimality for optimal control problems
  governed by maximal monotone operators.
\newblock {\em Journal of Differential Equations}, 260(1):733--757, 2016.

\bibitem{Brog03}
B.~Brogliato.
\newblock Some perspectives on the analysis and control of complementarity
  systems.
\newblock {\em {IEEE} Transactions on Automatic Control}, 48(6):918--935, 2003.

\bibitem{BrogGoel11}
B.~Brogliato and D.~Goeleven.
\newblock Well-posedness, stability and invariance results for a class of
  multivalued {L}ur'e dynamical systems.
\newblock {\em Nonlinear Analysis Series A: Theory, Methods \& Applications},
  74:195--212, 2011.

\bibitem{BrogThib10}
B.~Brogliato and L.~Thibault.
\newblock Existence and uniqueness of solutions for non-autonomous
  complementarity dynamical systems.
\newblock {\em Journal of Convex Analysis}, 17(3 \& 4):961--990, 2010.

\bibitem{Brow63}
F.E. Browder.
\newblock Variational boundary value problems for quasi-linear elliptic
  equations of arbitrary order.
\newblock {\em Proc. Natl. Acad. Sci., USA}, 50:31--37, 1963.

\bibitem{camlibel:02a}
M.K. Camlibel, W.P.M.H. Heemels, and J.M. Schumacher.
\newblock Consistency of a time-stepping method for a class of piecewise-linear
  networks.
\newblock {\em {IEEE Transactions on Circuits and Sys\-tems--I: Fundamental
  Theory and Applications}}, 49(3):349--357, 2002.

\bibitem{CamlHeem02}
M.K. Camlibel, W.P.M.H. Heemels, and J.M. Schumacher.
\newblock On linear passive complementarity systems.
\newblock {\em European Journal of Control}, 8(3):220--237, 2002.

\bibitem{camlibel:14}
M.K. Camlibel, L.~Iannelli, and F.~Vasca.
\newblock Passivity and complementarity.
\newblock {\em Mathematical Programming, Ser.~A}, 145:531--563, 2014.

\bibitem{j06}
M.K. Camlibel and J.M. Schumacher.
\newblock Existence and uniqueness of solutions for a class of piecewise linear
  dynamical systems.
\newblock {\em {Linear Algebra and its Applications}}, 351-352:147--184, 2002.

\bibitem{CamlSchu16}
M.K. Camlibel and J.M. Schumacher.
\newblock Linear passive systems and maximal monotone mappings.
\newblock {\em Mathematical Programming, Ser.~B}, 157(2):397--420, 2016.

\bibitem{Comb18}
P.L. Combettes.
\newblock Monotone operator theory in convex optimization.
\newblock {\em Mathematical Programming}, 170(1):177--206, 2018.

\bibitem{cottle:92}
R.W. Cottle, J.-S. Pang, and R.E. Stone.
\newblock {\em The Linear Complementarity Problem}.
\newblock Academic Press, Boston, 1992.

\bibitem{EdmoThib05}
J.F. Edmond and L.~Thibault.
\newblock Relaxation of an optimal control problem involving a perturbed
  sweeping process.
\newblock {\em Mathematical Programming, Ser.~B}, 104:347--373, 2005.

\bibitem{EdmoThib06}
J.F. Edmond and L.~Thibault.
\newblock {BV} solutions of nonconvex sweeping process differential inclusion
  with perturbation.
\newblock {\em J. Differential Equations}, 226:135--179, 2006.

\bibitem{pang}
F.~Facchinei and J.-S. Pang.
\newblock {\em Finite-Dimensional Variational Inequalities and Complementarity
  Problems I}.
\newblock {Springer}, New York, 2003.

\bibitem{j16}
L.~Han, A.~Tiwari, M.K. Camlibel, and J.-S. Pang.
\newblock Convergence of time-stepping schemes for passive and extended linear
  complementarity systems.
\newblock {\em SIAM Journal on Numerical Analysis}, 47(5):3768--3796, 2009.

\bibitem{Heem00}
W.P.M.H. Heemels, J.M. Schumacher, and S.~Weiland.
\newblock Linear complementarity systems.
\newblock {\em SIAM J. Appl. Math.}, 60:1234--1269, 2000.

\bibitem{HuPapa00}
S.~Hu and N.~Papageorgiu.
\newblock {\em Handbook of Multivalued Analysis}, volume II: Applications of
  {\em Mathematics and Its Applications}.
\newblock Springer, 2000.

\bibitem{JourVilc18}
A.~Jourani and E.~Vilches.
\newblock A differential equation approach to implicit sweeping processes.
\newblock {\em Journal of Differential Equations}, 2018.

\bibitem{kandilakis96}
D.A. Kandilakis.
\newblock Nonlinear evolution equations involving time-dependent
  subdifferentials of opposite sign.
\newblock {\em Nonlinear Analysis: Theory, Methods \& Applications},
  27(3):321--326, 1996.

\bibitem{kartsatos84}
A.G. Kartsatos and M.E. Parrott.
\newblock Functional evolution equations involving time dependent maximal
  monotone operators in banach spaces.
\newblock {\em Nonlinear Analysis: Theory, Methods \& Applications},
  8(7):817--833, 1984.

\bibitem{Kato67}
T.~Kato.
\newblock Nonlinear semigroups and evolution equations.
\newblock {\em J.Math.Soc.Japan}, 19(4):508--520, 1967.

\bibitem{Kacu60}
R.I. Ka\v{c}urovski\u{i}.
\newblock Monotone operators and convex functionals.
\newblock {\em Uspekhi Mat. Nauk}, 15:213--215, 1960.

\bibitem{KrejRoch11}
P.~Krejci and T.~Roche.
\newblock Lipschitz continuous data dependence of sweeping processes in {BV}
  spaces.
\newblock {\em Discrete and Continuous Dynamical Systems series B},
  15(3):637--350, 2011.

\bibitem{Kubo88}
M.~Kubo.
\newblock Subdifferential operator approach to nonlinear age-dependent
  population dynamics.
\newblock {\em Japan J. Appl. Math.}, 5:225--256, 1988.

\bibitem{KunzMont97}
M.~Kunze and M.D.P.~Monteiro Marques.
\newblock {BV} solutions to evolution problems with time-dependent domains.
\newblock {\em Set-Valued Analysis}, 5:57--72, 1997.

\bibitem{KunzMont00}
M.~Kunze and M.D.P.~Monteiro Marques.
\newblock An introduction to {M}oreau's sweeping process.
\newblock In B.Brogliato, editor, {\em Impacts in Mechanical Systems. Analysis
  and Modelling}, volume 551 of {\em Lecture Notes in Physics}, pages 1--60.
  Springer-Verlag, Berlin, 2000.

\bibitem{Kuttler00}
K.~Kuttler.
\newblock Evolution inclusions for time dependent families of subgradients.
\newblock {\em Applicable Analysis}, 76(3-4):185--201, 2000.

\bibitem{LaksLeel81}
V.~Lakshmikantham and S.~Leela.
\newblock {\em Nonlinear Differential Equations in Abstract Spaces}.
\newblock Pergamon Press, 1981.

\bibitem{LeraLion65}
J.~Leray and J.-L.Lions.
\newblock Quelques r\'esultats de {V}\v{i}sik sur les probl\`emes elliptiques
  nonlin\'eaires par les m\'ethodes de {Minty-Browder}.
\newblock {\em Bull. Soc. Math. France}, 93:97--107, 1965.

\bibitem{Mint62}
G.J. Minty.
\newblock Monotone (nonlinear) operators in {Hilbert} space.
\newblock {\em Duke Math. J.}, 29:341--346, 1962.

\bibitem{Mont93}
M.D.P. {Monteiro Marques}.
\newblock {\em Differential Inclusions in Nonsmooth Mechanical Problems: Shocks
  and Dry Friction}, volume~9 of {\em Progress in Nonlinear Differential
  Equations and their Applications}.
\newblock Birkh\"auser, 1993.

\bibitem{More77}
J.J. Moreau.
\newblock Evolution problem associated with a moving convex set in a {H}ilbert
  space.
\newblock {\em J. Differential Equations}, 26:347--374, 1977.

\bibitem{otani94}
M.~Otani.
\newblock Nonlinear evolution equations with time-dependent constraints.
\newblock {\em Advances in Mathematical Sciences and Applications}, 3, 01 1994.

\bibitem{pang2012}
J.-S. Pang, L.~Han, G.~Ramadurai, and S.~Ukkusuri.
\newblock A continuous-time linear complementarity system for dynamic user
  equilibria in single bottleneck traffic flows.
\newblock {\em Mathematical Programming, Ser.~A}, 133(1):437--460, 2012.

\bibitem{pang2009}
J.-S. Pang and D.~Stewart.
\newblock Solution dependence on initial conditions in differential variational
  inequalities.
\newblock {\em Mathematical Programming, Ser.~B}, 116(1):429--460, 2009.

\bibitem{PangStew08}
J.-S. Pang and D.E. Stewart.
\newblock Differential variational inequalities.
\newblock {\em Mathematical Programming, Ser.~A}, 113:345--424, 2008.

\bibitem{Papa90}
N.~Papageorgiou.
\newblock On evolution inclusions associated with time-dependent convex
  subdifferential.
\newblock {\em Commentationes Mathematicae Universitatis Carolinae},
  031(3):517--527, 1990.

\bibitem{PariBoyd13}
N.~Parikh and S.~Boyd.
\newblock Proximal algorithms.
\newblock {\em Foundations and Trends in Optimization}, 1(3):123--231, 2013.

\bibitem{Pave87}
N.H. Pavel.
\newblock {\em Nonlinear Evolution Operators and Semigroups: {A}pplications to
  Partial Differential Equations}, volume 1260 of {\em Lecture Notes in
  Mathematics}.
\newblock Springer, Berlin, 1987.

\bibitem{PeypSori10}
J.~Peypouquet and S.~Sorin.
\newblock Evolution equations for maximal monotone operators: Asymptotic
  analysis in continuous and discrete time.
\newblock {\em J. Convex Analysis}, 17(3 \& 4):1113--1163, 2010.

\bibitem{Phel93}
R.R. Phelps.
\newblock {\em Convex Functions, Monotone Operators and Differentiability}.
\newblock Springer-Verlag, Berlin, 2nd edition, 1993.

\bibitem{Recu15}
V.~Recupero.
\newblock {BV} continuous sweeping processes.
\newblock {\em Journal of Differential Equations}, 259:4253--4272, 2015.

\bibitem{Rock76}
R.T. Rockafellar.
\newblock Monotone operators and the proximal point algorithm.
\newblock {\em SIAM Journal on Control and Optimization}, 14(5):877--898, 1976.

\bibitem{Rockafellar98}
R.T. Rockafellar and J.-B. Wets.
\newblock {\em Variational Analysis}.
\newblock A Series of Comprehensive Studies in Mathematics 317. Springer, 1998.

\bibitem{rockafellar:70}
{R.T. Rockafellar}.
\newblock {\em Convex Analysis}.
\newblock Princeton University Press, Princeton, New Jersey, 1970.

\bibitem{rudin}
W.~Rudin.
\newblock {\em Principles of Mathematical Analysis}.
\newblock International Series in Pure and Applied Mathematics. McGraw-Hill
  Book Co., New York, third edition, 1976.

\bibitem{RyuBoyd16}
E.K. Ryu and S.~Boyd.
\newblock A primer on monotone operator methods: {S}urvey.
\newblock {\em Appl. Comput. Math.}, 15(1):3--43, 2016.

\bibitem{schumacher2004}
J.-M. Schumacher.
\newblock Complementarity systems in optimization.
\newblock {\em Mathematical Programming, Ser.~B}, 101(1):263--295, 2004.

\bibitem{Show93}
R.E. Showalter.
\newblock {\em Monotone Operators in Banach Space and Nonlinear Partial
  Differential Equations}.
\newblock American Mathematical Society, Providence, RI, 1997.

\bibitem{SiddManc02}
A.H. Siddiqi, P.~Manchanda, and M.~Brokate.
\newblock On some recent developments concerning {M}oreau's sweeping process.
\newblock In A.H. Siddiqi and M.~Kocvara, editors, {\em Proceedings of the 1st
  International Conference on Industrial and Applied Mathematics of the Indian
  Subcontinent: Trends in Industrial and Applied Mathematics}, Applied
  Optimization, pages 339--354. Kluwer Academic Publishers, 2002.

\bibitem{Simo08}
S.~Simons.
\newblock {\em From {Hahn-Banach} to Monotonicity}.
\newblock Springer-Verlag, Berlin, 2nd edition, 2008.

\bibitem{TanwBrog14}
A.~Tanwani, B.~Brogliato, and C.~Prieur.
\newblock Stability and observer design for {L}ur'e systems with multivalued,
  non-monotone, time-varying nonlinearities and state jumps.
\newblock {\em SIAM J. Control and Optimization}, 56(2):3639--3672, 2014.

\bibitem{TanwBrog18}
A.~Tanwani, B.~Brogliato, and C.~Prieur.
\newblock Well-posedness and output regulation for implicit time-varying
  evolution variational inequalities.
\newblock {\em SIAM J. Control and Optimization}, 56(2):751--781, 2018.

\bibitem{Tols17}
A.A. Tolstonogov.
\newblock Existence and relaxation of solutions for a subdifferential inclusion
  with unbounded perturbation.
\newblock {\em J. Math. Anal. Appl.}, 477(1):269--288, 2017.

\bibitem{Vlad91}
A.A. Vladimirov.
\newblock Nonstationary dissipative evolution equations in a {H}ilbert space.
\newblock {\em Nonlinear Analysis}, 17:499--518, 1991.

\bibitem{yamazaki05}
N.~Yamazaki.
\newblock Nonlinear evolution equations with time-dependent constraints.
\newblock {\em Hokkaido University Preprint Series in Mathematics}, 696:1--16,
  2005.

\bibitem{Zara60}
E.H. Zarantonello.
\newblock Solving functional equations by contractive averaging.
\newblock Technical Report Tech. report no. 160, University of Wisconsin,
  Madison, USA, 1960.

\bibitem{Zeid90}
E.~Zeidler.
\newblock {\em Nonlinear Functional Analysis and its Applications II/B --
  Nonlinear Monotone Operators}.
\newblock Springer-Verlag, NewYork, 1990.

\end{thebibliography}

\end{document}